\documentclass[DIV=14,letterpaper]{scrartcl}
\usepackage{tikz}
\usetikzlibrary{arrows}
\usetikzlibrary{calc}
\usetikzlibrary{matrix}
\usetikzlibrary{patterns}
\usepackage{amsmath,amssymb,amsfonts,amsthm,subfig} % Typical maths resource packages
\usepackage{verbatim}
\usepackage{extarrows}
\usepackage{mathrsfs}
\usepackage{mathtools}
\usepackage{color}                    % For creating coloured text and background
\usepackage{hyperref}                 % For creating hyperlinks in cross references
\usepackage{etoolbox}
\usepackage{graphicx}
\usepackage{tikz-cd}
\newtheorem{theorem}{Theorem}[section]
\newtheorem{corollary}{Corollary}[theorem]
\newtheorem{lemma}[theorem]{Lemma}
\newtheorem{remark}{Remark}

\def \D{\mathbb{D}}

\newcommand{\bey}{\begin{eqnarray}}
\newcommand{\eey}{\end{eqnarray}}

\newcommand{\beq}{\begin{equation}}
\newcommand{\eeq}{\end{equation}}
\theoremstyle{plain}% default

\theoremstyle{definition}

\theoremstyle{remark}
\newtheorem{exam}{\hspace{6mm}\textbf{Example}}[section]
\newtheorem{rem}{\hspace{6mm}{Remark}}[section]

\newcommand{\dif}{\textrm{d}}

\patchcmd{\subequations}{\alph}{.\arabic}{}{} % change display format of eq. counter

\vspace{5mm}

\begin{document}

\date{}
\title{Moving Mesh with Streamline Upwind Petrov-Galerkin (MM-SUPG) Method for Time-dependent Convection-Dominated Convection-Diffusion Problems}
\author{Xianping Li% 
\thanks{College of Integrative Sciences and Arts, Arizona State University, AZ, 85212, USA (\textit{xianping@asu.edu})}
\and Matthew McCoy%
\thanks{Valorem Reply, Kansas City, MO, 64157, USA (\textit{matthewmccoy6@gmail.com})}
}
\maketitle

\vspace{10pt}

\begin{abstract}
Time-dependent convection-dominated convection-diffusion problems are considered. We develop a moving mesh streamline upwind Petrov-Galerkin (MM-SUPG) method by combining residual-based SUPG stabilization with a metric-based moving mesh PDE (MMPDE) approach. 
The key feature of the method is the interaction between mesh adaptation and stabilization: the evolving mesh modifies both the direction and magnitude of the SUPG stabilization through the local element geometry, while stabilization improves robustness in convection-dominated regimes.
For isotropic diffusion, numerical results show that the proposed method reduces spurious oscillations and provides improved resolution of sharp layers compared with fixed mesh methods, while yielding accuracy comparable to moving mesh finite element methods without SUPG in some cases. 
For anisotropic diffusion, we introduce a weighted tensor that incorporates both the diffusion tensor and the convection field, and construct a metric tensor via intersection to guide mesh adaptation.
Under suitable assumptions on the elementwise tensors and mesh geometry, we establish sufficient conditions for the discrete maximum principle (DMP) of the fully discrete scheme. The analysis is based on quantitative bounds for the convection terms, 
and requires structural conditions on the alignment between the velocity field and the diffusion tensor. Numerical experiments demonstrate that the proposed metric improves monotonicity properties and reduces undershoots, while maintaining overall accuracy.

\end{abstract}

% \begin{keywords}
%     Convection-diffusion equation, SUPG, moving mesh, MMPDE, anisotropic diffusion, discrete maximum principle
% \end{keywords}

\vspace{10pt}

%%%%%% Description of the problem %%%%%%
\section{Introduction}
We consider the time-dependent convection-diffusion equation defined over the time-independent domain $\Omega\subset \mathbb{R}^2$ with a Lipschitz boundary and denote $\Omega_T = \Omega\times [0,T]$. The equation is given by:
\begin{subequations}
    \label{modeleqn}
    \begin{align}
        \label{modelequation1}
        &u_t- \varepsilon \nabla \cdot (\mathbb{D}\nabla u) + {\boldsymbol b}({\boldsymbol x, t})\cdot \nabla u=f({\boldsymbol x},t)\textrm{\hspace{0.38in} in \hspace{0.25in}} \Omega\times (0,T], \\
        \label{Dirichletboundary}
        &u({\boldsymbol x},t) = g({\boldsymbol x},t) \textrm{\hspace{1.8in} on \hspace{0.25in}} \partial \Omega\times [0,T], \\
        \label{modeleqnboundary}
        &u({\boldsymbol x},0) = u_0({\boldsymbol x}) \textrm{\hspace{1.85in} in\hspace{0.25in}} \Omega\times\{t = 0\}.
    \end{align}
\end{subequations}
Here, $\mathbb{D} = \mathbb{D}({\boldsymbol x})$ is the diffusivity tensor that is symmetric and positive-definite (SPD), 
${\boldsymbol b}$ represents an incompressible velocity field, while $g$, $u_0$, and $f$ are sufficiently smooth functions. 
The parameter $\varepsilon \in (0,1]$ is the scaling factor for the diffusivity coefficient. 
When $\D$ has the same eigenvalues, the problem is isotropic. In the anisotropic case, the eigenvalues of $\D$ are different, 
which indicates that diffusion in one direction is faster than in the other. 

Our focus lies in the convection-dominated scenario where $\varepsilon$ is much smaller than $\|{\boldsymbol b}\|_{L^{\infty}}$. 
Such problems arise in various fields, including pollutant dispersal in rivers, atmospheric pollution, the Stefan problem, turbulent transport, and are 
particularly intriguing in the investigation of the Navier--Stokes equations with large Reynolds numbers, see \cite{morton2019revival} for relevant applications.
When $\varepsilon$ is small, non-physical oscillations tend to arise in the solution unless an extremely fine spatial-temporal discretization is employed, incurring substantial computational costs. 
To enhance stability, artificial diffusion in the direction of ${\boldsymbol b}$ is proposed in \cite{brooks1982streamline} to be added to the standard Galerkin finite element method. 

The streamline upwind Petrov-Galerkin (SUPG) method has proven to be particular effective in numerically addressing boundary and interior layer problems in both time-dependent and time-independent scenarios \cite{di2018anisotropic,erath2019optimal,picasso2009adaptive,ganesan2017ale}.
The motivation behind SUPG is to employ upwinding of the test functions along the streamlines, aiming to enhance the diffusion of solution. This concept is most clearly illustrated in the 1D time-independent case using finite difference stencils: 
\begin{align}
    \frac{\partial u}{\partial x} = -\varepsilon \frac{\partial^2 u}{\partial x^2}.
\end{align}
When applying the finite difference method with only the center stencil, one obtains a numerical approximation that, depending on $\varepsilon$, tends to be under-diffusive. 
Conversely, utilizing only an upwind stencil leads to excessive diffusion, known as over-diffusion. As the exact solution typically lies somewhere in between these extremes, 
it prompts the adoption of upwinding the center scheme \cite{fries2004review}. This method was initially proposed by Brooks and Hughes for incompressible 
Navier-Stokes problems \cite{brooks1982streamline}, and has since proven to be fruitful in studying those equations \cite{franca1992stabilized,gelhard2005stabilized}. 
Comparison among different methods for the steady-state form of equation (\ref{modeleqn}) have been performed in \cite{augustin2011assessment}.

Efforts have been made to numerically solve equation (\ref{modelequation1}) for time-independent solutions, particularly focusing on steady flows \cite{di2018anisotropic, nguyen2009adaptive}, by employing adaptive mesh methods 
in conjunction with SUPG, especially when $\varepsilon \ll 1$. 
Notably, mesh adaptation has emerged as a potent tool in addressing convection-dominated problems. Our investigation targets a specific class of problems characterized by two key properties: 
artificial (non-physical) oscillations and sharp layers, encompassing both interior and boundary layers. Consequently, our interest lies in devising a mesh strategy grounded in anisotropic mesh adaptation techniques. 

Several notable methodologies have been explored for solving equation (\ref{modeleqn}). For instance, in \cite{bause2002uniform}, Lagrange-Galerkin error analysis demonstrated first-order uniform convergence in $\varepsilon$. 
Additionally, $\varepsilon$-uniform convergence theory for finite element method with homogeneous boundary conditions was investigated in \cite{bause2004uniform}. Carey et al. introduced adaptive mesh refinement (AMR) in \cite{carey2004some}, 
while Chinchapatnam et al. explored unsymmetric and symmetric meshless collocation techniques with radial basis functions in \cite{chinchapatnam2006unsymmetric}. 
Of particular interest are the so-called SOLD (Spurious Oscillations at Layers Diminishing) schemes, which have been extensively studied in \cite{john2008finite,john2009finite,john2007spurious}. These schemes involve the addition of either 
isotropic or anisotropic artificial diffusion to the Galerkin formulation. Furthermore, in \cite{cheng2020local}, local discontinuous Galerkin (LDG) methods were employed on a class of layer-adapted Shishkin and Bakhvalov-type 
meshes, accompanied by error analysis, to approximate equation (\ref{modeleqn}). 

For time-dependent problems, Wang et al. utilized the singular boundary method with $f =0$ in \cite{wang2017singular}, employing an exponential variable transform along with the Fourier transform. 
Additionally, unfitted finite element methods have gained traction for tackling equation (\ref{modeleqn}) in time-dependent domains \cite{von2020unfitted}. These types of unfitted FEMs are known by various names such as XFEM, CutFEM, Finite Cell Method, and TraceFEM (refer to the 
sources in \cite{von2020unfitted} along with \cite{hansbo2016cut, zahedi2017space}). It is worth noting, as mentioned in \cite{badia2018aggregated}, that unfitted techniques may encounter conditioning difficulties, which are addressed 
using specialized finite element spaces and aggregation technique. 

Recent developments include variational multiscale stabilized formulations combined with adaptive mesh strategies (e.g., \cite{zhang2024mmvmfem}), as well as residual-based stabilization on anisotropic meshes. 
Compared with these approaches, the present work focuses on integrating SUPG stabilization with a metric-based moving mesh PDE framework, with particular emphasis on discrete maximum principle (DMP) preservation for anisotropic diffusion.

In all the above mentioned methods, only isotropic diffusivity tensors were considered. However, in many practical applications, the diffusivity tensor is anisotropic.
When the mesh is not aligned in the principal diffusion direction, non-physical oscillations may occur when using standard computational methods \cite{li2010anisotropic}. 
Special care must be taken to ensure that the numerical solution satisfies the discrete maximum principle (DMP).
For problem with only diffusion, metric tensors has been developed in \cite{li2010anisotropic} for time-independent problems to adapt the mesh such as $\mathbb{M}_{DMP} = \theta_K \mathbb{D}_K^{-1}$ based on DMP satisfaction 
and $\mathbb{M}_{DMP+adap} = C_K (det(\mathbb{D}))^{\frac{1}{d}} \mathbb{D}_K^{-1}$ based on DMP satisfaction as well as minimization of interpolation errors.
For time-dependent problems, condition for time step size is also developed in \cite{li2013maximum} to ensure that the numerical solutions satisfy DMP.
Lu et al. investigated the satisfaction of DMP for anisotropic diffusion-convection-advection problems in \cite{lu2014maximum}.
However, it is even more challenge for numerical solutions of convection-dominated convection-diffusion problems to satisfy the DMP, which will be addressed in this paper.

In this paper, we propose to integrate the moving mesh partial differential equation (MMPDE) method with streamline upwind Petro-Galerkin (SUPG) method and call it as moving mesh SUPG (MM-SUPG) method. It works for both time-independent and time-dependent flows 
characterized by initial sharp layers, and can deal with both isotropic and anisotropic diffusivity tensors. 
It has been observed that, even for isotropic diffusivity, oscillations can persist with adapted meshes \cite{linss2001numerical}. The challenge in numerically solving equation (\ref{modeleqn}) becomes evident particularly when $\varepsilon \ll 1$ and 
the initial profile experiences convection over time or when the profile features sharp regular, parabolic, corner, or interior layers that develop or exist initially.
Meanwhile, when $\D$ is anisotropic, non-physical solutions may appear in the numerical solution, which is a violation of DMP \cite{li2010anisotropic,li2013maximum}. 
To address this issue, we introduce a metric tensor to control the mesh movement and develop conditions for mesh and time step size so that the numerical solution satisfies DMP.

The novelty of the proposed MM-SUPG method lies not merely in combining moving mesh adaptation with SUPG stabilization, but in analyzing and exploiting their interaction. 
The effectiveness of SUPG stabilization depends strongly on mesh geometry, as both the streamline derivative ${\boldsymbol b} \cdot \nabla u_h$ and the stabilization parameter $\tau_K$ are mesh-dependent. 
Meanwhile, the moving mesh method dynamically aligns elements according to solution features and anisotropic diffusion.
This coupling leads to a nontrivial interplay: mesh alignment influences the direction and strength of SUPG stabilization, which in turn affects oscillation control and monotonicity. 
In particular, for anisotropic diffusion, the proposed metric tensor construction balances alignment with both the convection field and diffusion tensor, enabling improved discrete maximum principle behavior. 
This interaction between mesh adaptation and stabilization is a central contribution of the present work.

The remainder of this paper is organized as follows. In Section 2, we give a brief introduction for SUPG and MMPDE, and propose our MM-SUPG method. Section 3 presents numerical examples for isotropic diffusivity case to demonstrate the effectiveness of our approach, while also 
delineating the constraints of the moving mesh method. In Section 4, we investigate the satisfaction of discrete maximum principle for anisotropic diffusion. Finally, a summary and some conclusions are presented in Section 5. 

\section{Moving Mesh with Streamline Upwind Petrov-Galerkin (MM-SUPG) Method}
Firstly, we give a brief introduction to SUPG method and MMPDE method. Then we explain in details our MM-SUPG method, the integration of MMPDE with SUPG. 

\subsection{SUPG method}
Consider the convection-diffusion equation (\ref{modelequation1}) with homogeneous Dirichlet boundary conditions (\ref{Dirichletboundary}) and initial condition (\ref{modeleqnboundary}). The SUPG method is a stabilized finite element method that introduces artificial diffusion to the Galerkin formulation.
Let $L^2(\Omega)$ and $H^k(\Omega)$, $k \ge 1$ denote the standard Lebesgue and Sobolev spaces, respectively. Let 
$H_g^1(\Omega) = \{v \in H^1(\Omega), \, v|_{\partial \Omega}=g\}$ and $H_0^1(\Omega) = \{v \in H^1(\Omega), \, v|_{\partial \Omega}=0\}$. The weak formulation of (\ref{modelequation1}) is to find $u\in H_g^1(\Omega)$ such that
\begin{equation}
    \label{weakform}
    \int_{\Omega} u_t v - \varepsilon \int_{\Omega} \nabla \cdot (\mathbb{D} \nabla u) v + \int_{\Omega} ({\boldsymbol b} \cdot \nabla u) v = \int_{\Omega} f v, \quad \forall v \in H_0^1(\Omega).
\end{equation}
We consider an affine family of simplicial triangulations $\{\mathcal{T}_h\}$ for $\Omega$ and denote the piecewise linear finite element space as $V^h = (v \in H^1(\Omega), \, v|_K \in \mathcal{P}_1^K, \forall K \in \mathcal{T}_h)$, where $\mathcal{P}_1^K$ is the linear polynomial space in element $K$.
Denote by $U^h = V^h \cap H_g^1(\Omega)$ and $U_0^h = V^h \cap H_0^1(\Omega)$. The standard Galerkin method is to find $u_h \in U^h$ such that
\begin{equation}
    \label{galerkin}
    \int_{\Omega} (u_h)_t v_h - \varepsilon \sum_{K \in \mathcal{T}_h} \nabla \cdot (\mathbb{D}_K \nabla u_h) v_h + \int_{\Omega}  ({\boldsymbol b} \cdot \nabla u_h) v_h = \int_{\Omega} f v_h, \quad \forall v_h \in U_0^h,
\end{equation}
where $\mathbb{D}_K$ is the restriction of $\mathbb{D}$ to element $K$. 

The essential idea of SUPG involves altering the test function space to operate within the realm of functions that upwind the trial functions along the streamlines. This can be represented as the space:
\begin{align}
    U_h^{\textrm{SUPG}}=\left\{v_h + \tau_K{\boldsymbol b}\cdot \nabla v_h \mid v_h\in U_0^h\right\},
\end{align}
where $\tau_K$ is the stabilization parameter.
Therefore, one may perceive the space $U_h^{\textrm{SUPG}}$ as a suitable stabilized test space for singularly perturbed problems. 
The weak formulation for SUPG method is given as follows: find $u_h \in U^h(\Omega)$ such that
\begin{align}
    \label{SUPGweakform}
    \nonumber
    \int_{\Omega} (u_h)_t v_h & - \varepsilon \sum_{K \in \mathcal{T}_h} \nabla \cdot (\mathbb{D}_K \nabla u_h) (v_h + \tau_K {\boldsymbol b} \cdot \nabla v_h) \\
    & + \sum_{K \in \mathcal{T}_h} ({\boldsymbol b} \cdot \nabla u_h) (v_h + \tau_K {\boldsymbol b} \cdot \nabla v_h) = \int_{\Omega} f (v_h + \tau_K {\boldsymbol b} \cdot \nabla v_h), \quad \forall v_h \in U_0^h.
\end{align}

% The standard choice of $\tau$ is given by
% \begin{equation}
%     \label{tau}
%     \tau = \frac{h}{2\|{\boldsymbol b}\|_{\infty}},
% \end{equation}

The stabilization parameter $\tau_K$ is chosen following classical SUPG formulations, where it scales with the local element size and inverse flow magnitude to balance numerical diffusion and consistency; 
see \cite{brooks1982streamline, franca1992stabilized, codina2000stabilization}, and related developments in moving mesh contexts \cite{di2018anisotropic}.
A conventional choice of $\tau_K$ on mesh $\mathcal{T}_h$ is 
\begin{align}
    \tau_{\mathcal{T}_h}^K =\begin{cases}
    c_0\textrm{diam}(K)_{\mathcal{T}_h}\ \ \textrm{if}\ \ \textrm{Pe}_{K}>1 \text{ convection-dominated}\\
    c_1\textrm{diam}(K)^2_{\mathcal{T}_h}\ \ \textrm{if}\ \ \textrm{Pe}_{K}\leq 1 \text{ diffusion-dominated}
    \end{cases},
\end{align}
where $\textrm{Pe}_K$ denotes the element-wise P\'eclet number given by 
\begin{align}
    \textrm{Pe}_K = \frac{\|{\boldsymbol b}\|_{\infty}\textrm{diam}(K)_{\mathcal{T}_h}}{2\varepsilon}.
\end{align}
$\tau_K$ determines local convection-dominated or diffusion-dominated regimes on the current mesh $\mathcal{T}_K$. 
There is flexibility in defining $\textrm{diam}(K)_{\mathcal{T}_h}$ and in selecting $c_0$ and $c_1$. This paper focuses on convection-dominated regimes characterized by large element P\'eclet numbers $P_{e_{K}} \gg 1$, 
which may arise due to small $\varepsilon$ and/or insufficient mesh resolution relative to the flow magnitude.. Typically, individual bilinear forms in the SUPG formulation are estimated, and $\tau_{\mathcal{T}_h}^K$ is derived based on those errors. 
In \cite{di2018anisotropic}, such a method was employed using error estimates in the discrete norm for time independent ${\boldsymbol b}$
\begin{align}
    \|w\|_{\varepsilon, k}^2 = \varepsilon \|\nabla w\|^2 +\sum_{K\in\mathcal{T}_K}\tau_{\mathcal{T}_h}^K\|{\boldsymbol b}\cdot \nabla w\|_{L^2(K)}^2, \, \forall w\in H_0^1(\Omega).
\end{align}
However, our method relies on different error estimates for defining the metric tensor along with a mesh functional approach. 
When a reaction/source term is present, stabilization parameters must be modified to account for its contribution to the residual, as discussed in \cite{codina2000stabilized,franca2000improved,knopp2002stabilized,lube2006residual}.

In the case of anisotropic mesh adaptation, $\textrm{diam}(K)_{\mathcal{T}_h}$ will be updated with each new mesh as the elements stretch and shrink. In this paper, we define $\tau_{\mathcal{T}_h}^K$ as 
\begin{align}
    \tau_{\mathcal{T}_h}^K &= \frac{\textrm{diam}(K)_{\mathcal{T}_h}}{2\|{\boldsymbol b}\|_{\infty}}\xi(\textrm{Pe}_K), 
    \label{taudef}
\end{align}
which is based on the stabilization parameter selected in \cite{nguyen2009adaptive}. Common choices of $\textrm{diam}(K)_{\mathcal{T}_h}$ include the length of the longest edge of the element $K$ projected onto ${\boldsymbol b}$ 
and the diameter of $K$ in the direction of ${\boldsymbol b}$. In this paper, for each $K\in \mathcal{T}_h$, we define  
\begin{align}
    \textrm{diam}(K) = \sup_{{\boldsymbol x}_1,{\boldsymbol x}_2\in K}\|{\boldsymbol x}_1-{\boldsymbol x}_2\|,
    \label{diameter}
\end{align}
and 
\begin{align}
    \xi(\textrm{Pe}_K) = \min\left\{1, \frac{\textrm{Pe}_K}{3}\right\}.
    \label{xifactor}
\end{align}

\subsection{Moving mesh partial differential equations (MMPDE)}

Moving mesh methods are instrumental in addressing domains where the solution exhibits excessively large gradient, such as 
sharp layers. Traditional methods like fitted operator methods, finite difference methods, finite element methods with uniform
or Shishkin-type meshes may not be able to catch the sharp gradient without a hyper-refined mesh or a priori knowledge about the solution. 
The MMPDE approach automatically concentrates mesh elements near the region where the solution exhibits large gradient, 
hence improves the computational efficiency. In this subsection, we give a brief overview of MMPDE. More details can be found in \cite{Huang_book}. 

We consider an affine family of simplicial triangulations $\{\mathcal{T}_h\}$ for $\Omega$. Denote by $(K, \mathcal{P}_1^K,\Sigma_K)$ the finite element space. 
Let $N$, $N_v$, and $N_t$ denote the number of mesh elements, the number of vertices, and the number of time steps, respectively. 
Let $\mathbb{M}=\mathbb{M}({\boldsymbol x})$ be a symmetric and uniformly positive definite metric tensor defined on $\Omega$ with 
$c_1\mathbb{I}\leq \mathbb{M}({\boldsymbol x})\leq c_2\mathbb{I}$ for all ${\boldsymbol x}\in \Omega$. 
We take the $\mathbb{M}$-uniform mesh approach that takes a nonuniform mesh and views it as a uniform one in the metric specified by $\mathbb{M}$. 
Essentially, $\mathbb{M}$ controls the concentration and alignment of the mesh elements.  

We consider two affine-equivalent finite elements, one for the reference elements and another for the physical mesh element. 
There exists an invertible affine mapping between the finite elements, denoted by $F_K$:
    \begin{align}
        (\hat{K},\mathcal{P}_1^{\hat{K}},\Sigma_{\hat{K}})	\xrightarrow{\quad F_K \quad}  (K,\mathcal{P}_1^K,\Sigma_K). 
    \end{align}
For each element $K\in\mathcal{T}_h$, denote by $F_K:\hat{K}\to K$ the affine mapping between $K$ and the reference element $\hat{K}$ so that $\lvert\hat{K}\rvert=1$.
Denote by ${\boldsymbol x}_j^K$ for $j=0$, 1, and 2 the vertices of $K$ and ${\boldsymbol \xi}_j^{K}$ the corresponding vertices of $\hat{K}$. Then the affine map $F_K$ maps ${\boldsymbol \xi}_j^{K}$ to ${\boldsymbol x}_j^K$, i.e.,  $F_K({\boldsymbol \xi}_j^{K})={\boldsymbol x}_j^K$ for $j=0$, 1, and 2. $\hat{K}$ is called the reference element and belongs to $\hat{\mathcal{T}}_{h}$, the reference mesh. 

Denote by $\mathbb{M}_K$ the volume-average of $\mathbb{M}$ over element $K$, i.e., 
    $$\mathbb{M}_K=\frac{1}{\lvert K\rvert}\int_K\mathbb{M}({\boldsymbol x})\dif{\boldsymbol x},$$
and
    $$\sigma_h=\sum_{K\in \mathcal{T}_h}\int_K\sqrt{\det(\mathbb{M}_K)}\dif {\boldsymbol x}=\sum_{K\in \mathcal{T}_h}\lvert K\rvert\sqrt{\det(\mathbb{M}_K)}.$$
$F_K$ and $\mathbb{M}$ give rise to the so-called equidistribution and alignment conditions in $\mathbb{R}^2$:
\begin{enumerate}
    \item equidistribution: $\lvert K\rvert \sqrt{\det(\mathbb{M}_K)}=\displaystyle\frac{\sigma_h}{N}$ for all $K\in\mathcal{T}_h$, and 
    \item alignment: $\displaystyle\frac{1}{2}\textrm{tr}\left(F_K'\mathbb{M}_K(F_K')^T\right)=\sqrt{{\textrm{det}}\left(F_K'\mathbb{M}_K(F_K')^T\right)}$ for all $K\in\mathcal{T}_h$, 
\end{enumerate}
where $F_K'$ is the Jacobian matrix of $F_K$. These two completely characterize a non-uniform mesh: the equidistribution condition controls the size and shape of 
mesh elements, while the alignment condition controls the orientation of mesh elements. The new physical mesh $\mathcal{T}_h$ is found by minimizing the energy functional 
(also called the mesh functional) $I$ given in discrete form by 
    \begin{align}
        I(\mathcal{T}_h)= \sum_{K\in\mathcal{T}_h}\lvert K\rvert G(\mathbb{J}_K,\det(\mathbb{J}_K),\mathbb{M}_K),
        \label{mmpde-functional}
    \end{align}
where $\mathbb{J}_K=(F_K')^{-1}$ and 
    \begin{align}
        \label{fun-G}
        G(\mathbb{J}_K,\det(\mathbb{J}_K),\mathbb{M}_K)=&\alpha\sqrt{\det(\mathbb{M}_K)}\left(\textrm{tr}\left(\mathbb{J}_K\mathbb{M}_K\mathbb{J}_K^T\right)\right)^p\\
        \nonumber
        &+(1-2\alpha)2^p\sqrt{\det(\mathbb{M}_K)}\left(\frac{\det(\mathbb{J}_K)}{\sqrt{\det(\mathbb{M}_K)}}\right)^p,
    \end{align}
with $\alpha\in (0,\frac{1}{2}]$ and $p>1$ as dimensionaless parameters. From numerical experiments, $\alpha = \frac{1}{3}$ and $p = \frac{3}{2}$ work well for most problems \cite{kolasinski}.

The moving mesh PDE (MMPDE) method is developed by Huang et al. \cite{Huang_book} to find a new physical mesh by minimizing $I(\mathcal{T}_h)$, which is equivalent to the steepest descent method \cite{zhang}. 
The MMPDE is a partial differential equation (also called the mesh PDE) that is given by 
    \begin{align}
        \frac{\dif {\boldsymbol x}_i}{\dif t} = -\frac{P_i}{\gamma}\left(\frac{\partial I}{\partial {\boldsymbol x}_i}\right)^T, \ i = 1, \dots, N_v, 
        \label{mmpde}
    \end{align}
where $\gamma$ is a user parameter meant to adjust the time scale of the mesh movement and $P_j = \det \left(\mathbb{M}(x_j)\right)^{\frac{p-1}{2}}$ is chosen so that the MMPDE is invariant 
under the scaling transformation of $\mathbb{M}$. 

The physical PDE and mesh PDE can be solved either simultaneously or alternatively. When the PDEs are solved simultaneously, the mesh 
can respond quickly to any change in the physical solution. However, this procedure is highly nonlinear. In our work, we employ the alternative approach. 
The mesh PDE is firstly solved based on the physical solution and the mesh at the current time step. Then mesh nodes are relocated based on the solutions from mesh PDE. 
Lastly, the physical PDE is solved on the new mesh for the next time step. The alternate procedure allows 
for flexibility and potential efficiency at each time level. One disadvantage of this approach is a lag in time, which can be addressed by employing multiple iterations on the same time level before advancing to the next. 
It is always beneficial to start with a initially adapted mesh, especially when sharp layers exist from the beginning.  

In this paper, for moving mesh for isotropic diffusion case, we use the metric tensor defined by 
    \begin{align}
        \mathbb{M}=\det\left(\mathbb{I}+\lvert H_K(u)\rvert\right)^{-\frac{1}{6}}\left(\mathbb{I}+\lvert H_K(u)\rvert\right),
        \label{metrictensordef}
    \end{align}
where $\mathbb{I}$ is the identity matrix, $H_K(u)$ is the recovered Hessian using least squares fitting to the values of the solution function at the mesh vertices, 
and $\lvert H_K(u)\rvert=\Lambda\textrm{diag}(\lvert\lambda_1\rvert,\lvert\lambda_2\rvert)\Lambda^T$ with $\Lambda\textrm{diag}(\lambda_1,\lambda_2)\Lambda^T$ being the eigen-decomposition of $H_K(u)$. The Hessian is 
computed $H(u)({\boldsymbol x}_j) \approx H(q)({\boldsymbol x}_j)$, where $q$ is the algebraic polynomial determined by least squares fitting. For completeness, some details for the least squares fitting 
method (see Section 5.3 in \cite{Huang_book}) are provided here. 

Specifically, if $u_j$ is the approximation of $u$ at node ${\boldsymbol x}_j=(x_j,y_j)$ on mesh $\mathcal{T}_h$. Let ${\boldsymbol x}_{j_i}$ for $i = 1,\dots,N_j$ be the $N_j$ neighboring vertices near ${\boldsymbol x}_j$ (including ${\boldsymbol x}_j$). 
Denote by ${\boldsymbol x}_j^c$ the center of these points and define  
\begin{align}
    H_j^x = \max_{i=1,\dots,N_j}\lvert x_{j_i}-x_j^c\rvert,\ \ H_j^y = \max_{i=1,\dots,N_j}\lvert y_{j_i}-y_j^c\rvert.
\end{align}
Then the polynomial $q$ is in (Legendre) polynomial form
\begin{align}
    q(x,y) = \sum_{k=0}^2\sum_{m =0}^{2-k} a_{k,m}P_k\left(\frac{x-x_j^c}{H_j^x}\right)P_m\left(\frac{y-y_j^c}{H_j^y}\right),
\end{align}
and is determined by
\begin{align}
    \min_{a_{0,0},\dots,, a_{2,0}}\sum_{i = 1}^{N_j}(q(x_{j_i},y_{j_i})-u_{j_i})^2. 
\end{align}
Then it is clear that $\nabla u({\boldsymbol x}_j) \approx \nabla q({\boldsymbol x}_j)$. This metric tensor is optimal in the $L^2$ norm of linear interpolation error on triangular meshes \cite{zhang}. 
Moreover, it is clearly symmetric and one may replace $\lvert H_K(u)\rvert$ with $\alpha^{-1}\vert H_K(u)\rvert$, where $\alpha$ is called the intensity parameter and helps ensure that $\mathbb{M}$ is positive definite. 
This monitor function is derived on anisotropic error estimates (see \cite{Huang_book}, chapter 5, for a more detailed discussion and the derivation of (\ref{metrictensordef}) and choice of $\alpha$).

\subsection{Moving mesh with Streamline Upwind Petrov-Galerkin (MM-SUPG) method}

Now we apply the MMPDE method to solve the SUPG formulation \eqref{SUPGweakform}. 
Write ${\boldsymbol u}_h({\boldsymbol x},t) = \sum\limits_{j =1}^{N_v} u_j(t)\varphi_j({\boldsymbol x})$, where $\{\varphi_j\}_{j=1}^{N_v}$ is the set of linear basis functions in $U_0^h$. 
Also, since $\varphi_j$ is a linear function, the element-wise integral average can be defined explicitly
\begin{align}
    \label{basisformula}
    \frac{1}{|K|}\int_K\varphi_j\varphi_i\dif{\boldsymbol x} = \begin{cases}
    \displaystyle \frac{1}{(d+1)(d+2)}\textrm{\ \ for \ \  }i\neq j\\
    \displaystyle \frac{2}{(d+1)(d+2)}\textrm{\ \ for \ \ }i=j
    \end{cases}. 
\end{align}

Approximate the global integral as 
\begin{align}
    \int_{\Omega}g({\boldsymbol x},t)\dif {\boldsymbol x}&= \sum_{K\in\mathcal{T}_h}|K|\int_Kg({\boldsymbol x},t)\dif {\boldsymbol x}\\
    &\approx\sum_{K\in\mathcal{T}_h}|K|\sum_{j=1}^m\lambda_jg({\boldsymbol \omega}_j),
\end{align}
where ${\boldsymbol \omega_j}$ are the quadrature nodes for element $K$ and $\lambda_j$ are weights such that $\sum_{j=1}^m\lambda_j=1$. 
For any tensor $\mathbb{F}$, we then define the element-wise average by 
\begin{align}
    \label{avg_tensor}
    \mathbb{F}_K = \frac{1}{|K|}\int_K \mathbb{F}({\boldsymbol x})\dif {\boldsymbol x} \approx \sum_{j=1}^N\lambda_j\mathbb{F}({\boldsymbol \omega_j}). 
\end{align}

The semi-discrete matrix formulation for (\ref{SUPGweakform}) is given by 
\begin{align}
    {\boldsymbol M}{\boldsymbol u}_h'(t)+{\boldsymbol A}(t){\boldsymbol u}_h(t)={\boldsymbol f}(t),
    \label{TimeODE}
\end{align}
where ${\boldsymbol u} = [u_1,\dots,u_{N_v}]^T$, ${\boldsymbol M}$, ${\boldsymbol A}$, and ${\boldsymbol f}$ are given by 
\begin{align}
    \label{MmatrixFormula}
    M_{ij} &= \sum_{K\in\mathcal{T}_h}(\varphi_j, \varphi_i)_{0,K} + \tau_K(\varphi_j,{\boldsymbol b}\cdot \nabla \varphi_i)_{0,K}, \\
    \label{Amatrixint}
    A_{ij}(t) &= \sum_{K\in\mathcal{T}_h}({\boldsymbol b}\cdot \nabla \varphi_j,\varphi_i)_{0,K} + |K|(\nabla \varphi_i)^T \varepsilon \mathbb{D}_K\nabla\varphi_j\\
    \nonumber 
     & \quad\quad +\tau_K\left(({\boldsymbol b}\cdot \nabla \varphi_j, {\boldsymbol b}\cdot \nabla \varphi_i)_{0,K} + (\varepsilon \mathbb{D}_K \nabla\varphi_j,\nabla ({\boldsymbol b}\cdot\nabla\varphi_i))_{0,K}\right), \\
    \label{Frightside}
     f_i(t) & = \sum_{K\in\mathcal{T}_h}(f, \varphi_i)_{0,K}+\tau_K(f,{\boldsymbol b}\cdot \nabla \varphi_i)_{0,K}.
\end{align}

We use a $\theta$-scheme to discretize in time. For (\ref{TimeODE}), we approximate
\begin{align}
    {\boldsymbol u}_h'(t_m) \approx \frac{{\boldsymbol u}_h(t_{m+1})-{\boldsymbol u}_h(t_m)}{\Delta t},\ \ m = 1,\dots, N_t.
\end{align}
Then the corresponding $\theta$-scheme is given by 
\begin{align}
    M\frac{{\boldsymbol u}_h(t_{m+1})-{\boldsymbol u}_h(t_m)}{\Delta t}&+\theta A{\boldsymbol u}_h(t_{m+1})+(1-\theta)A{\boldsymbol u}_h(t_m)\label{discretizedtime}\\
    \nonumber 
    & = \theta {\boldsymbol f}(t_{m+1})+(1-\theta){\boldsymbol f}(t_m),\ \ m = 1,\dots, N_t.
\end{align}
Grouping the ${\boldsymbol u}_h(t_{m+1})$ terms, we may rewrite (\ref{discretizedtime}) as 
\begin{align}
    \left(\frac{M}{\Delta t}+\theta A\right){\boldsymbol u}_h(t_{m+1}) = \theta {\boldsymbol f}(t_{m+1})+(1-\theta){\boldsymbol f}(t_m)+\left(\frac{M}{\Delta t}-(1-\theta)A\right){\boldsymbol u}_h(t_m),
\end{align}
for $m = 1,\dots, N_t$. Hence, we may write a system of ODEs 
\begin{align}
    \label{ODEsystem}
    \tilde{{\boldsymbol A}}{\boldsymbol u}_h(t_{m+1})=\tilde{\boldsymbol f}^{m+1},\ \ m = 1,\dots, N_t,
\end{align}
where
\begin{align}
        \label{Amatrix}
        \tilde{{\boldsymbol A}}&={\boldsymbol M}+\theta \Delta t{\boldsymbol A}, \\
        \label{fvector}
        \tilde{\boldsymbol f}^{m+1}&=\theta\Delta t {\boldsymbol f}(t_{m+1})+\Delta t(1-\theta){\boldsymbol f}(t_m)+\left({\boldsymbol M}-\Delta t(1-\theta){\boldsymbol A}\right){\boldsymbol u}_h(t_m).
\end{align}
In this work, we use $\theta = 0.5$, which corresponds to the Crank-Nicolson method. Fixed time step is used in the computations. In Section 4 we derive conditions on $\Delta t$ so that DMP is satisfied. 

To generate the adapted mesh at each physical time level, we minimize the discrete meshing functional $I(\mathcal{T}_h)$ defined in \eqref{mmpde-functional} by solving the mesh PDE \eqref{mmpde}. 
The derivatives of the meshing functional are assembled element-by-element following the MMPDE framework in \cite{Huang_book,zhang} and the mesh PDE can be written in the form
\begin{align}
    \frac{\dif {\boldsymbol x}_i}{\dif t} = \frac{P_i}{\gamma} \sum_{K \in \omega_i} |K| {\boldsymbol v}_{i_K}^K , \ i = 1, \dots, N_v, 
    \label{mmpde-discretized}
\end{align}
where $\omega_i$ is the element patch associated with the $i$-th vertex, $i_K$ is the local index of vertex $i$ in element $K$, and ${\boldsymbol v}_{i_K}^K$ denotes the local mesh velocity contributed by element $K$.
Here the local velocities are obtained from the derivatives of $G(\mathbb{J}_K,\det(\mathbb{J}_K),\mathbb{M}_K)$ defined in \eqref{fun-G} with respect to $\mathbb{J}_K$ and $\det(\mathbb{J}_K)$, where $\mathbb{J}_K=(F_K')^{-1}$. 
Boundary vertices are constrained to remain on $\partial \Omega$, thus boundary nodes are allowed to slide along the boundary with corner nodes being fixed.

To obtain the new mesh, we integrate the discretized mesh PDE \eqref{mmpde-discretized} from the current mesh $\mathcal{T}_h^m$ to the updated mesh $\mathcal{T}_h^{m+1}$.
In our implementation, this system of ODEs is solved in pseudo-time using MATLAB \textit{ode15s} solver. After the integration is completed, the new mesh $\mathcal{T}_h^{m+1}$ is obtained directly from the updated vertex coordinates.

A key point is that the matrices ${\boldsymbol M}$ and ${\boldsymbol A}$, and the stabilization parameter $\tau_K$ are all assembled on the current physical mesh. 
Since the mesh changes in time, these quantities are recomputed after each mesh update. In particular, $\tau_K$ is evaluated on the updated mesh, so the amount of streamline stabilization adapts dynamically to the local mesh geometry.
We emphasize that we do not solve the MMPDE and the physical PDE fully simultaneously. Instead, we use a staggered procedure: \\
    \indent (1) fix the physical solution and update the mesh; \\
    \indent (2) freeze the new mesh and solve the SUPG system for the physical solution. \\
This decouples the nonlinear mesh movement from the stabilized finite element solution and is considerably simpler and more robust in practice. The algorithm is given as follows. \\

\noindent\rule{\textwidth}{0.4pt}
\vspace{1em}
\textbf{Algorithm 1: MM-SUPG Method} \\
    for $m=1, ..., N_t$, do: \\
    \indent $t_{m+1} = t_m + \Delta t$, 
    \begin{itemize}
        \item [(1)] Metric construction. \\
            Using ${\boldsymbol u}_h^m$ on $\mathcal{T}_h^m$, construct the metric tensor $\mathbb{M}^m = \mathbb{M}({\boldsymbol u}_h^n)$. 
        \item [(2)] Mesh movement. \\
            Fixing $\mathbb{M}^m$, solve the MMPDE for the mesh coordinates from ${\boldsymbol x}^m$ to an updated vertex set ${\boldsymbol x}^{m+1}$, yielding the new mesh $\mathcal{T}_h^{m+1}$.
        \item [(3)] Solution interpolation. \\
            Interpolate ${\boldsymbol u}_h^m$ from $\mathcal{T}_h^m$ onto $\mathcal{T}_h^{m+1}$, denoted by $\tilde{{\boldsymbol u}}_h^m$. 
        \item [(4)] Recompute geometric quantities. \\
            Using $\tilde{{\boldsymbol u}}_h^m$ on $\mathcal{T}_h^{m+1}$, compute $\textrm{diam}(K)_{\mathcal{T}_h}$, $\tau_{\mathcal{T}_h}^K$ and basis-function gradients.
        \item [(5)] Assemble SUPG system. \\
            Assemble the mesh-dependent matrices ${\boldsymbol M}^{M+1}$, ${\boldsymbol A}^{m+1}$, and ${\boldsymbol f}^{m+1}$.
        \item [(6)] Solve for physical solution. \\
            Solve the following linear system for the solution at $t_{m+1}$.
            $${\boldsymbol M}^{m+1}+\theta \Delta t{\boldsymbol A}^{m+1} {\boldsymbol u}_h^{m+1} = \theta\Delta t {\boldsymbol f}^{m+1}+\Delta t(1-\theta){\boldsymbol f}^{m}+\left({\boldsymbol M}^{m+1}-\Delta t(1-\theta){\boldsymbol A}^{m+1} \right){\boldsymbol u}_h^m.$$
    \end{itemize}

\noindent\rule{\textwidth}{0.4pt}
\vspace{1em}

The effectiveness of MM-SUPG stems from the complementary roles of its components: moving mesh adaptation reduces interpolation error by aligning elements with solution features, 
while SUPG stabilization suppresses spurious oscillations along streamlines. Together, they yield improved stability without excessive numerical diffusion.

\section{Results for Isotropic Diffusion} 
Now, we present some numerical examples for isotropic diffusion case where we take the diffusivity tensor $\mathbb{D}$ as an identity matrix. 
The value of diffusivity is given by $\varepsilon$. We consider two different examples and compare results from four different methods: 
fixed mesh finite element method (FM-FEM) without SUPG, fixed mesh finite element method with SUPG (FM-SUPG), moving mesh finite element method (MM-FEM) without SUPG (which is the MMPDE method), 
and our moving mesh with SUPG (MM-SUPG) method. We denote the MMPDE method as MM-FEM for consistency in what follows. 

Note that for $\varepsilon = O(1)$, all methods produce similar results, confirming that stabilization and mesh adaptation are primarily beneficial in convection-dominated regimes.

% We note that computational time is larger for implementation of SUPG, especially when using moving mesh. 

%%%%%%%%%%%%%%%%%%%%%%          Example 1       %%%%%%%%%%%%%%%%%%%%%%%%%%%%%%%%%%%%%%%%
\begin{exam}
\label{example1}
This example is similar to the simulated example presented in \cite{xie2019error}. The domain is $\Omega = (0,10)\times(-1.25, 1.25)$ and an initial smooth solution is given by
\begin{align}
    u_0(x,y) = \textrm{e}^{-\frac{(x-1)^2+y^2}{0.25}}+\textrm{e}^{-\frac{(x-1)^2+(y-0.5)^2}{0.25}}+\textrm{e}^{-\frac{(x-1)^2+(y+0.5)^2}{0.25}}.
\end{align}
The velocity field is given by 
\begin{align}
    {\boldsymbol b} = \left(1-\exp(\theta x)\cos(2\pi y), \displaystyle\frac{\theta}{2\pi}\exp(\theta x)\sin(2\pi y)\right),
    \label{example1flow}
\end{align}
with 
\begin{align}
    \theta = \frac{1}{2\varepsilon}-\sqrt{\left(\frac{1}{2\epsilon}\right)^2+4\pi^2}.
    \label{example1theta}
\end{align}

\begin{figure}[h]
    \centering
    \hbox{
    \begin{minipage}[t]{3.0in}
    \includegraphics[width=3.0in]{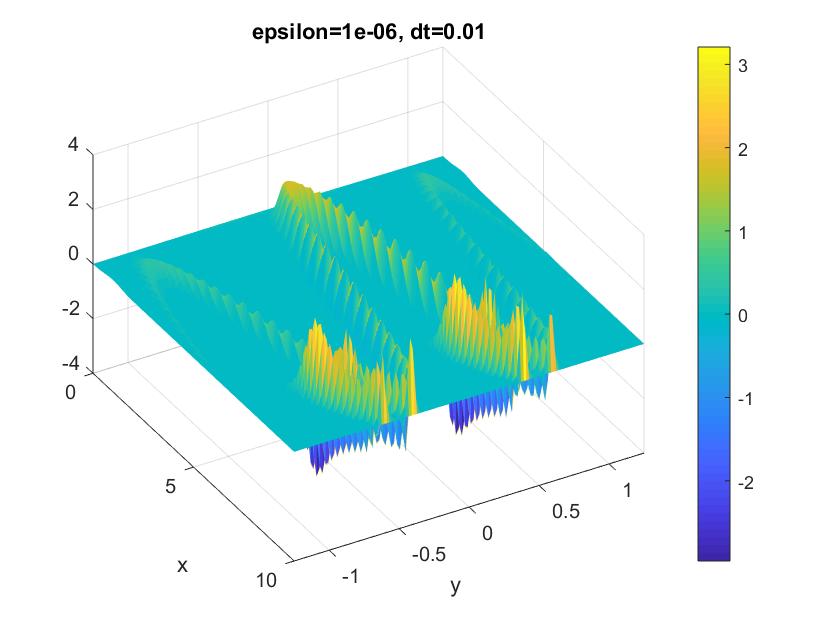}
    \centerline{(a): FM-FEM}
    \end{minipage}
    \hspace{10mm}
    \begin{minipage}[t]{3.0in}
    \includegraphics[width=3.0in]{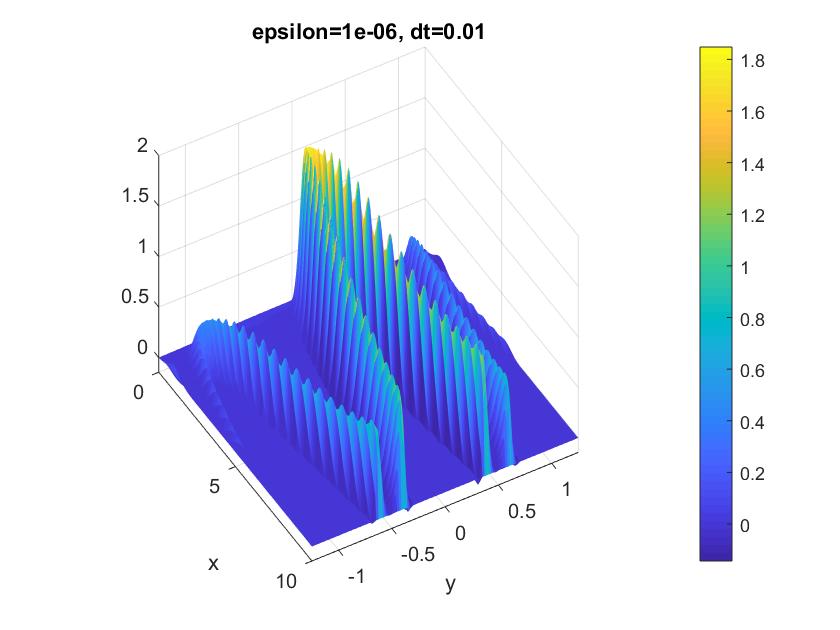}
    \centerline{(b): FM-SUPG}
    \end{minipage}
    }
    \hbox{
    \begin{minipage}[t]{3.0in}
    \includegraphics[width=3.0in]{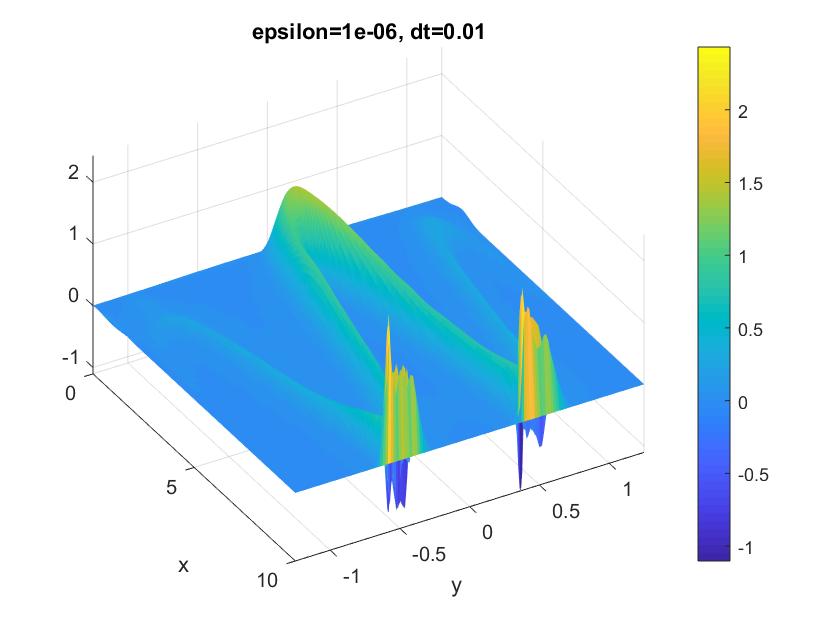}
    \centerline{(c): MM-FEM}
    \end{minipage}
    \hspace{10mm}
    \begin{minipage}[t]{3.0in}
    \includegraphics[width=3.0in]{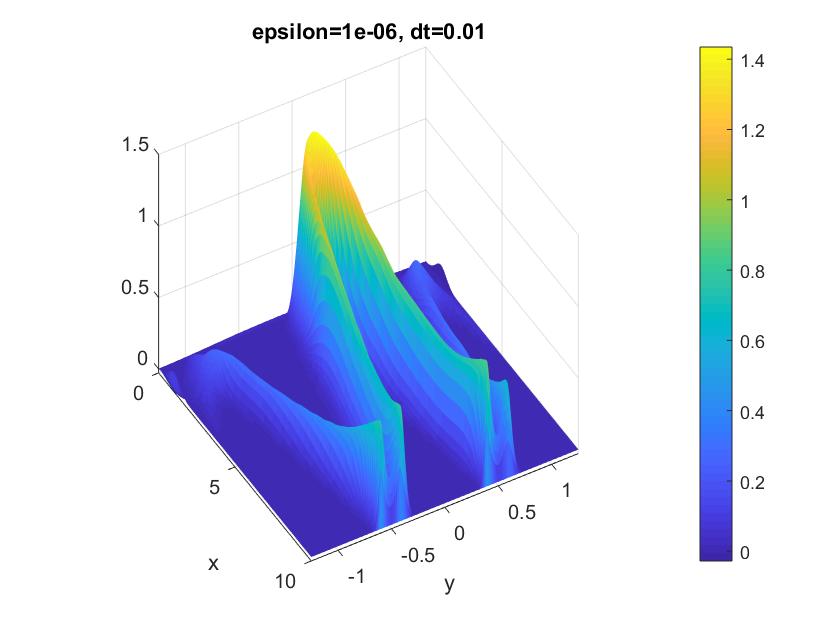}
    \centerline{(d): MM-SUPG}
    \end{minipage}
    }
    \caption{Example \ref{example1} -- side view of four different methods with $\Delta t = 0.01$, ${\boldsymbol b}$ is given by (\ref{example1flow}), $N = 32,768$ at $T=5.0$.}
    \label{ex1-sideview}
\end{figure}

For this example, $\varepsilon = 10^{-6}$ is used in the computation, which gives $\theta \approx -3.948\times 10^{-5}$. 
We compute the diffusion of the initial solution from 0 to 5 seconds with Dirichlet boundary conditions. Time step $\Delta t=10^{-2}$ and $N = 32,768$ elements are used in the computations.
Allowing the convection process to be carried out, we see fine layers form on the interior and boundary of $\Omega$. 

The results are shown in Figure \ref{ex1-sideview} with side view profiles.
From Figures \ref{ex1-sideview}(a) and (c), we see that for the classic Galerkin method, both FM-FEM and MM-FEM, the numerical solutions have non-physical oscillations form near the outflow boundary at $x = 10$. 
Figure \ref{ex1-sideview}(b) shows that FM-SUPG resolves the boundary layers present in FM-FEM and MM-FEM but has artificial oscillations on the interior of $\Omega$.    
On the other hand, Figure \ref{ex1-sideview}(d) shows that MM-SUPG results in the smoothest surface among all four methods, and both interior oscillations and oscillations near the outflow have been resolved. 
The computations are also performed with smaller $\Delta t$ as $5 \times 10^{-3}$, and the results are shown in Figure \ref{ex1-smaller-dt}. 
As can be seen, even with smaller time step, MM-FEM still has layers at the outflow while MM-SUPG performs much better.

\begin{figure}[thb]
    \centering
    \hbox{
    \begin{minipage}[t]{3.0in}
    \includegraphics[width=3.0in]{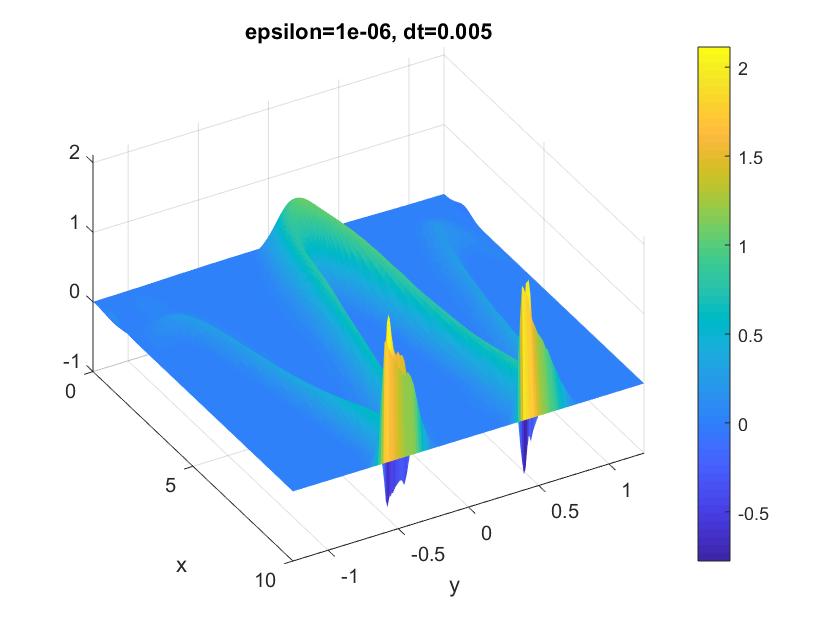}
    \centerline{(a): MM-FEM}
    \end{minipage}
    \hspace{10mm}
    \begin{minipage}[t]{3.0in}
    \includegraphics[width=3.0in]{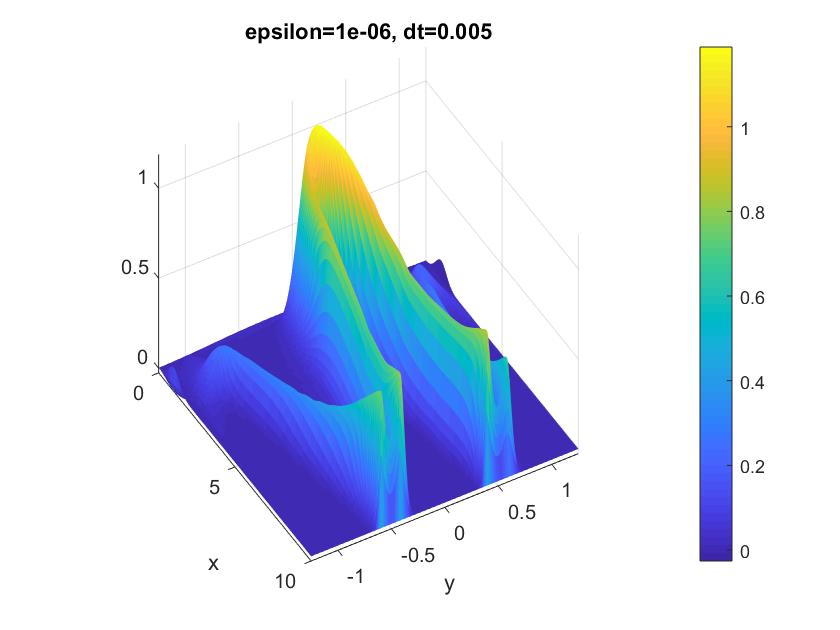}
    \centerline{(b): MM-SUPG}
    \end{minipage}
    }
    \caption{Example \ref{example1} -- side view of MM-FEM and MM-SUPG methods with $\Delta t = 0.005$, ${\boldsymbol b}$ is given by (\ref{example1flow}), $N = 32,768$ at $T=5.0$.}
    \label{ex1-smaller-dt}
\end{figure}

\end{exam}

%%%%%%%%%%%%%%%%%%%%%%          Example 2       %%%%%%%%%%%%%%%%%%%%%%%%%%%%%%%%%%%%%%%%

\begin{exam}
\label{example2}
    Here we consider another example (taken from \cite{john2008finite}) with exact solution 
    \begin{align}
        \nonumber
        u(x,y,t) = 1&6\sin(\pi t)x(1-x)y(1-y)\\
        &\times \left(\frac{1}{2}+\frac{1}{\pi}\arctan\left(2\varepsilon^{-1/2}(0.25^2-(x-0.5)^2-(y-0.5)^2)\right)\right),
    \end{align}
    with Dirichlet boundary conditions. We compute to time $T=0.5$ with time step $\Delta t=10^{-3}$. 
    The severity of the interior layer is proportional to $\varepsilon$ and its thickness is $\sqrt{\varepsilon}$. 
    We also use ${\boldsymbol b}=(2,3)$ and $\varepsilon = 10^{-6}$. 

    For coarse mesh with fewer elements (that is, smaller $N$), we observe artificial oscillations in the direction of ${\boldsymbol b}$, including MM-FEM, for example, with $N = 8192$. 
    Increasing the number of mesh elements helps to improve the results for MM-FEM. Figure \ref{ex2-topview} shows the top view of the solutions at $T=0.5$ for the four different methods with $N = 32768$, 
    while Figure \ref{ex2-sideview} displays the side view of the solution profiles. As can be seen, significant smearing occurs in the solution obtained using the FM-FEM and FM-SUPG, though the latter shows an improvement. 
    On the other hand, both the MM-FEM and MM-SUPG perform well with $N = 32768$ elements in the mesh. This is supported by the $H^1$-seminorm errors in Table \ref{ex2-timeindependent-H1} and numerical solutions shown in Figure \ref{ex2-sideview}. 
    The $L^2$-norm errors in Table \ref{ex2-timeindependent-L2} shows that MM-SUPG performs comparable even better than MM-FEM for this example.
    We observe pronounced mesh movement in the direction of ${\boldsymbol b}$ for MM-FEM and MM-SUPG, as shown in Figure \ref{ex2-mesh}. 
    The rate of convergence in the $H^1$-seminorm and $L^2$ norm, based on errors from Tables \ref{ex2-timeindependent-H1} and \ref{ex2-timeindependent-L2} are shown in Figure \ref{ex2-timeindependent}.
    The reference slopes correspond to expected convergence rates: $O(h)$ in $H^1$ and $O(h^{3/2}$ in $L^2$ using MMPDE framework for convection-dominated problems.
    
    \begin{figure}[thb]
        \centering
        \hbox{
        \begin{minipage}[t]{3.0in}
        \includegraphics[width=3.0in]{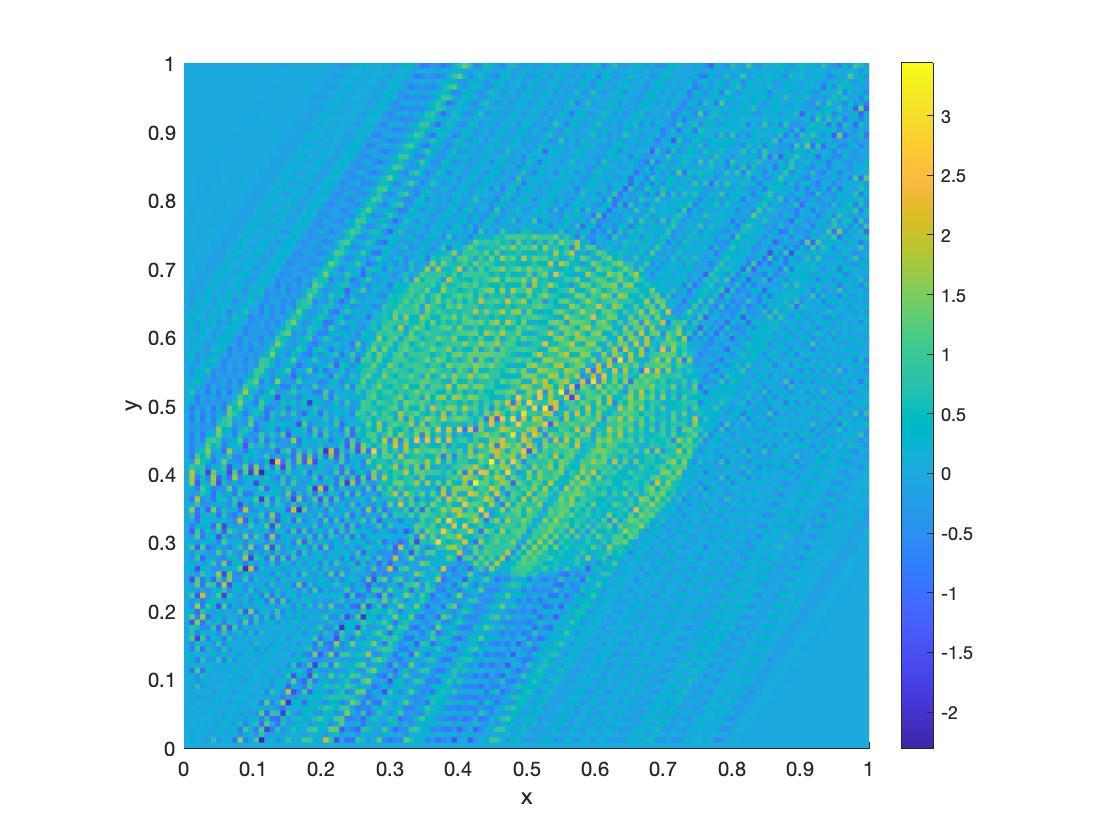}
        \centerline{(a): FM-FEM}
        \end{minipage}
        \hspace{10mm}
        \begin{minipage}[t]{3.0in}
        \includegraphics[width=3.0in]{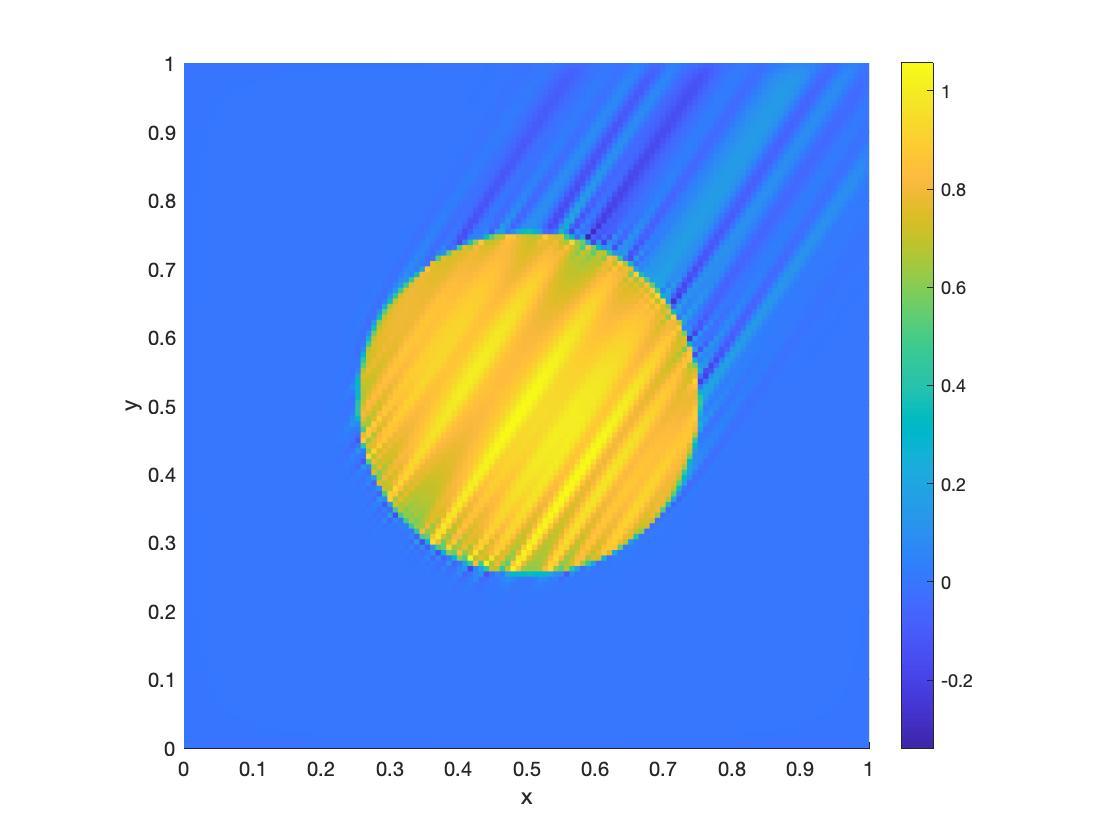}
        \centerline{(b): FM-SUPG}
        \end{minipage}
        }
        \hbox{
        \begin{minipage}[t]{3.0in}
        \includegraphics[width=3.0in]{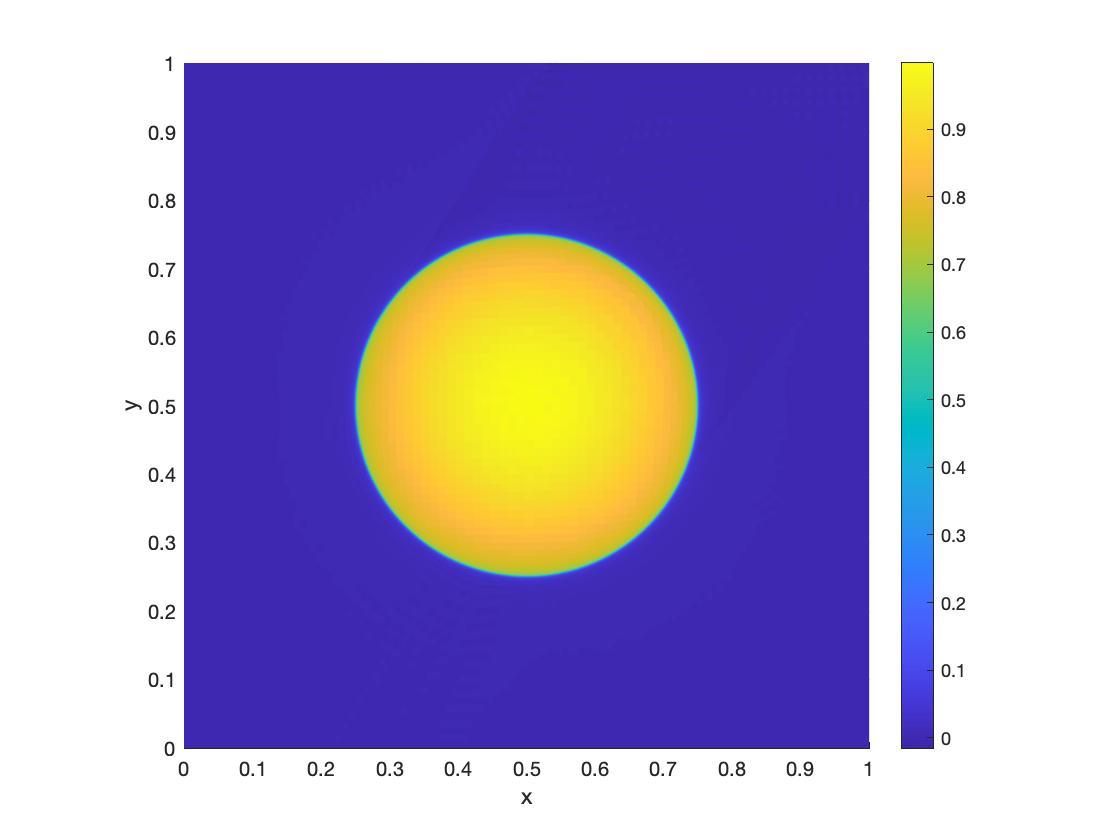}
        \centerline{(c): MM-FEM}
        \end{minipage}
        \hspace{10mm}
        \begin{minipage}[t]{3.0in}
        \includegraphics[width=3.0in]{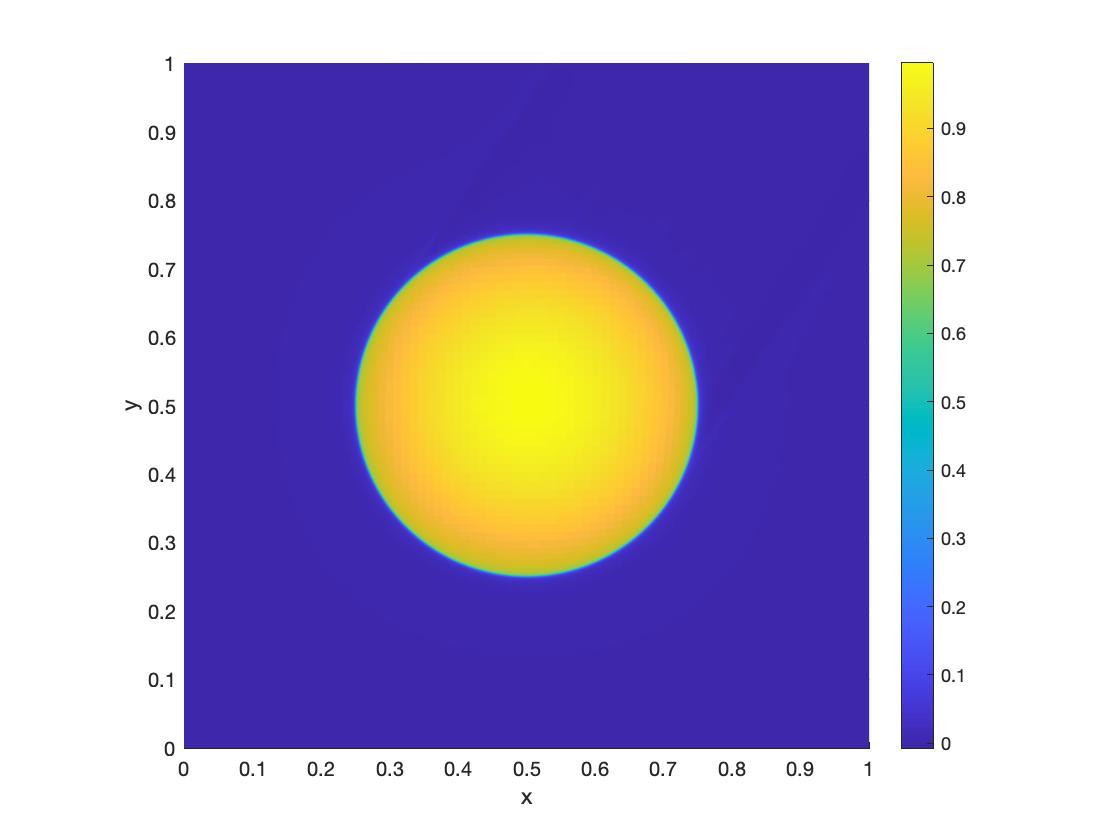}
        \centerline{(d): MM-SUPG}
        \end{minipage}
        }
        \caption{Example \ref{example2} -- top view with four different methods with $\Delta t = 0.001$ for ${\boldsymbol b} = (2,3)$, $N = 32768$ at $T=0.5$. 
        The observed discreteness in (a) is due to the use of piecewise linear finite elements on a fixed mesh where numerical oscillation is significant, 
        which results into obvious color changes in the top view of the solution.}
        \label{ex2-topview}
    \end{figure}

    \begin{figure}[thb]
        \centering
        \hbox{
        \begin{minipage}[t]{3.0in}
        \includegraphics[width=3.0in]{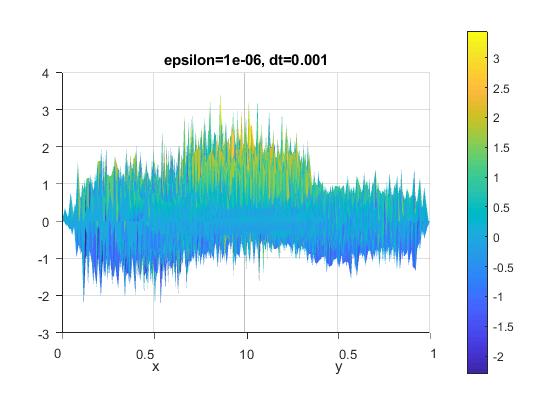}
        \centerline{(a): FM-FEM}
        \end{minipage}
        \hspace{10mm}
        \begin{minipage}[t]{3.0in}
        \includegraphics[width=3.0in]{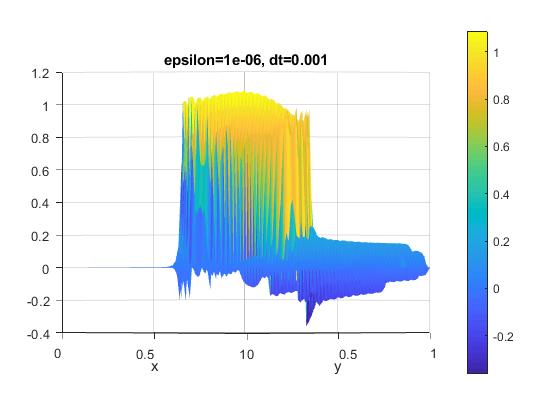}
        \centerline{(b): FM-SUPG}
        \end{minipage}
        }
        \hbox{
        \begin{minipage}[t]{3.0in}
        \includegraphics[width=3.0in]{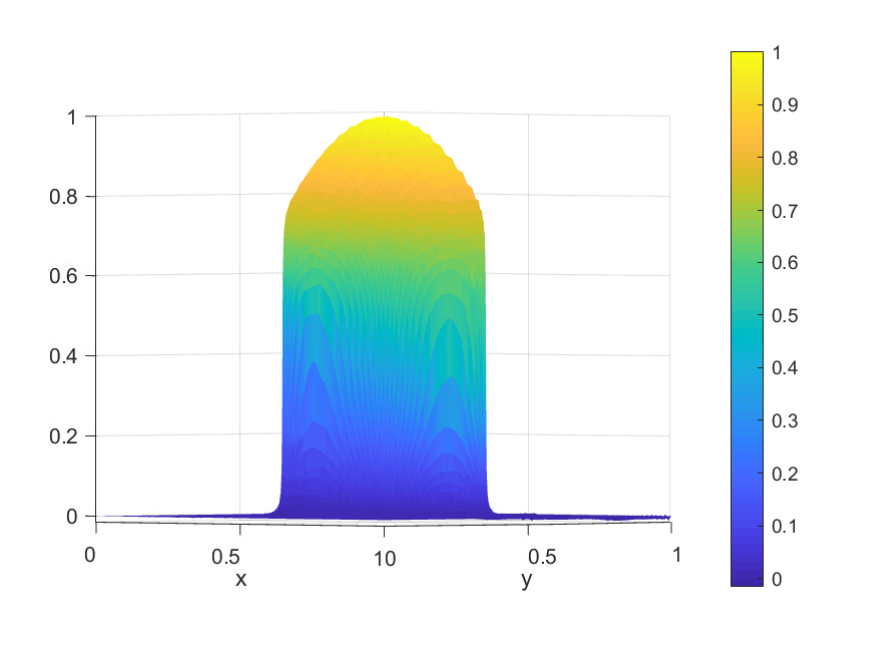}
        \centerline{(c): MM-FEM}
        \end{minipage}
        \hspace{10mm}
        \begin{minipage}[t]{3.0in}
        \includegraphics[width=3.0in]{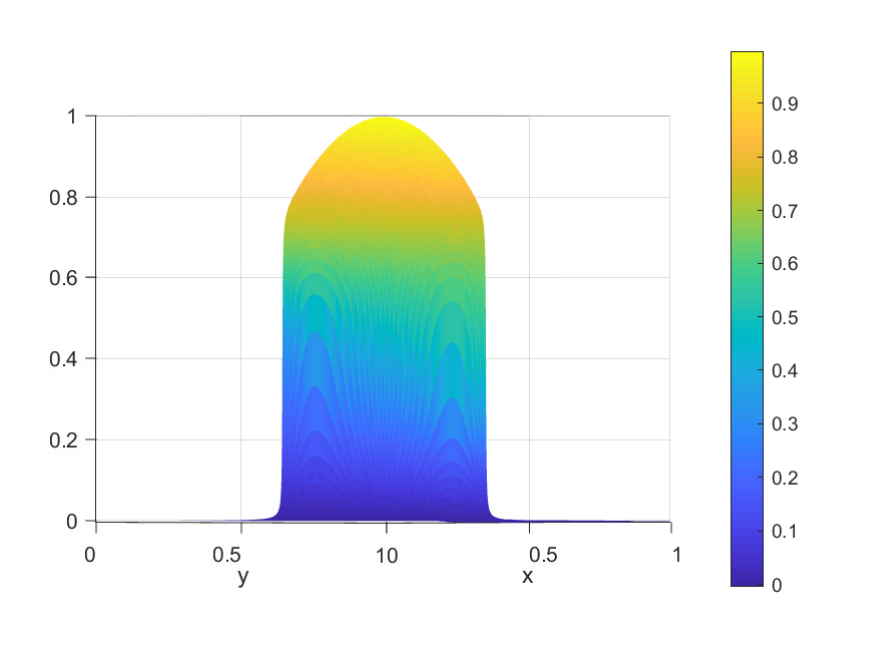}
        \centerline{(d): MM-SUPG}
        \end{minipage}
        }
        \caption{Example \ref{example2} -- side view with the four different methods with $\Delta t = 0.001$, ${\boldsymbol b} = (2,3)$, $N = 32768$ at $T=0.5$.}
        \label{ex2-sideview}
    \end{figure}
      
    \begin{figure}[thb]
        \centering
        \hbox{
        \begin{minipage}[t]{3.0in}
        \includegraphics[width=3.0in]{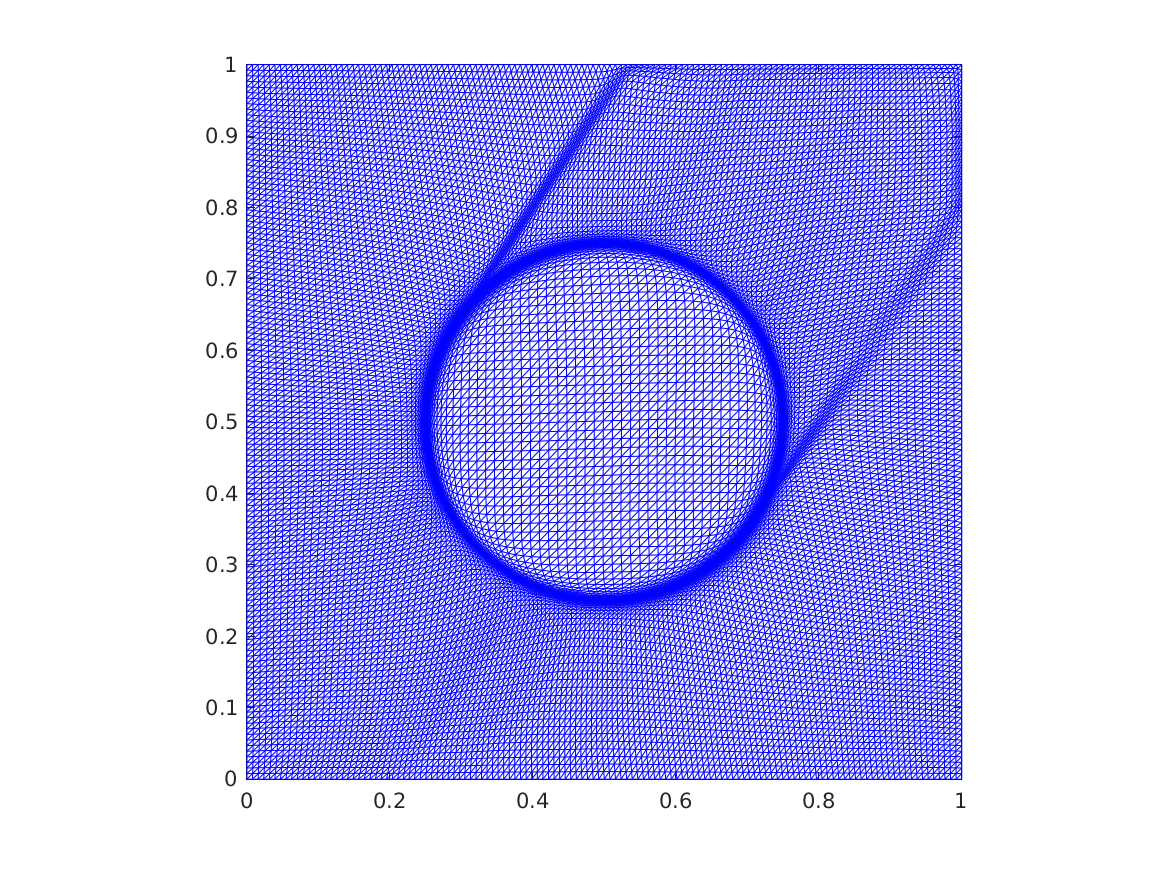}
        \centerline{(a): MM-FEM}
        \end{minipage}
        \hspace{10mm}
        \begin{minipage}[t]{3.0in}
        \includegraphics[width=3.0in]{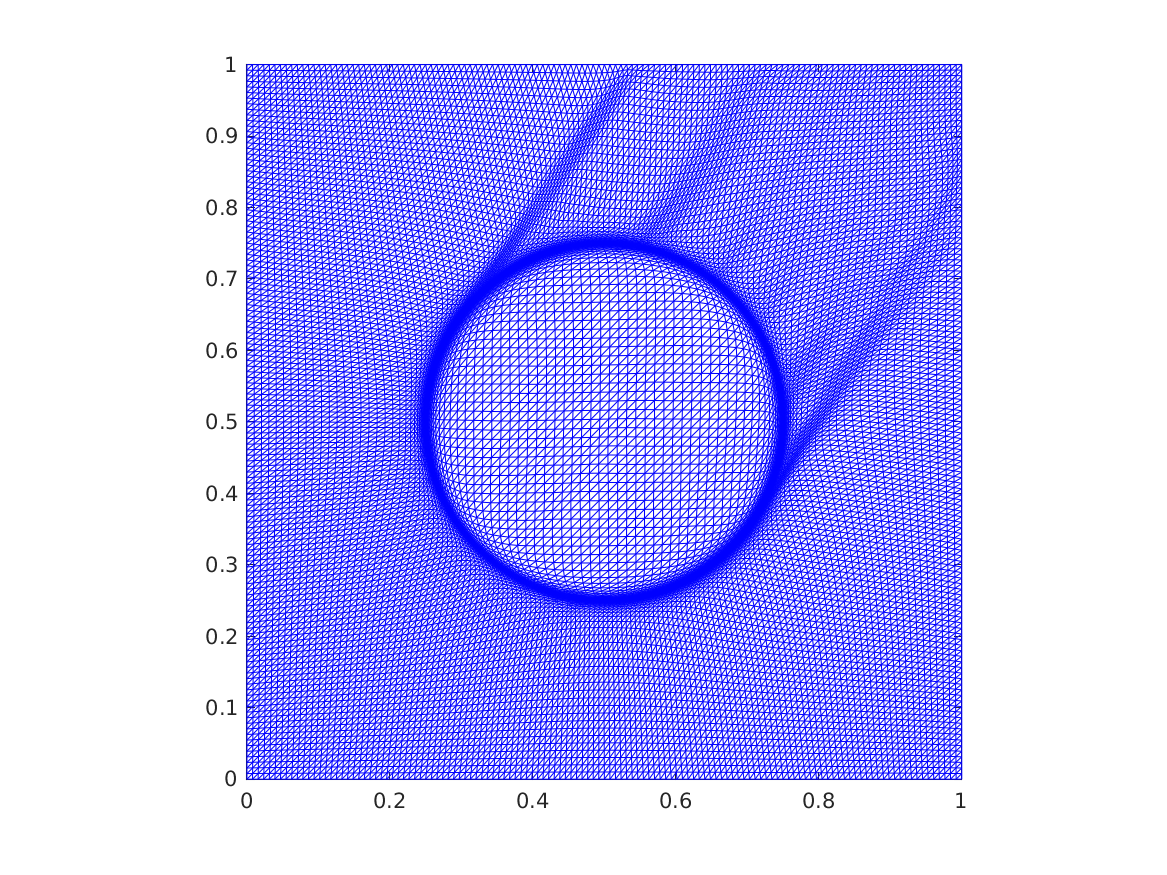}
        \centerline{(b): MM-SUPG}
        \end{minipage}
        }
        \caption{Example \ref{example2} -- Mesh movement for $T=0.5$ and ${\boldsymbol b} = (2,3)$, $N=32768$ at $T=0.5$.}
        \label{ex2-mesh}
    \end{figure}

    \begin{figure}[thb]
        \centering
        \hbox{
        \begin{minipage}[t]{3.0in}
        \includegraphics[width=3.0in]{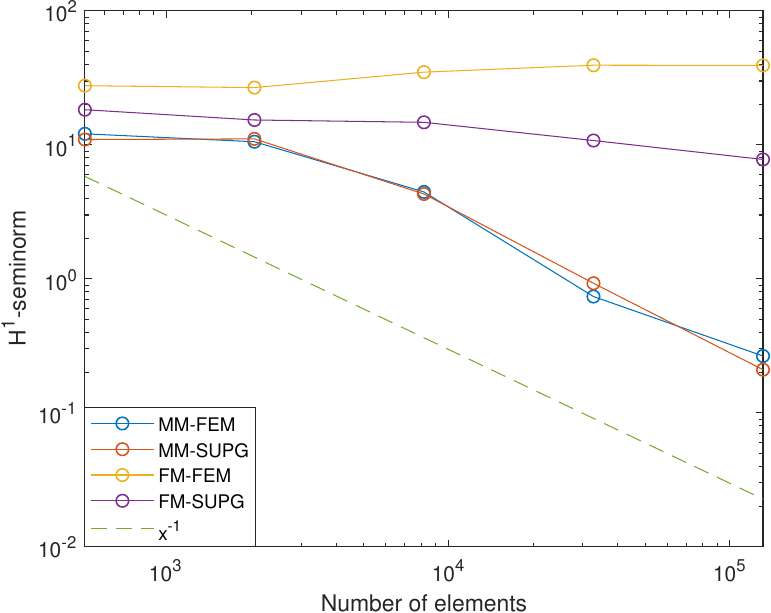}
        \centerline{(a): $H^1$-seminorm error.}
        \end{minipage}
        \hspace{10mm}
        \begin{minipage}[t]{3.0in}
        \includegraphics[width=3.0in]{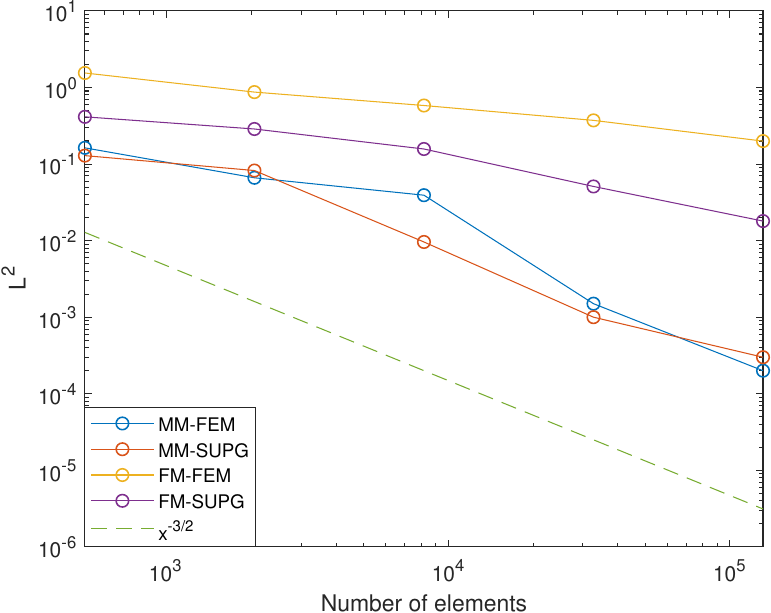}
        \centerline{(b): $L^2$ norm error.}
        \end{minipage}
        }
        \caption{Example \ref{example2} -- Errors corresponds to Tables \ref{ex2-timeindependent-H1} and \ref{ex2-timeindependent-L2}. $H^1$-seminorm and $L^2$ norm errors with $\varepsilon = 10^{-6}$ and $\Delta t = 0.001$ with flow field ${\boldsymbol b} = (2,3)$ at $T=0.5$.} 
        \label{ex2-timeindependent}
    \end{figure}
    
    \begin{table}[h!]
        \centering
        \caption{Example \ref{example2} -- $H^1$-seminorm errors for $\varepsilon=10^{-6}$ at $\Delta t = 0.001$ for ${\boldsymbol b}=(2,3)$.}
        \label{ex2-timeindependent-H1}
        \begin{tabular}{|l|l|l|l|l|}
        \hline
        $N$    & MM-FEM  & MM-SUPG & FM-FEM  & FM-SUPG \\ \hline
        512    & 12.1038 & 10.9715 & 27.6628 & 18.3147 \\ \hline
        2048   & 10.5624 & 11.0912 & 26.7879 & 15.3485 \\ \hline
        8192   & 4.4580  & 4.3027  & 34.9201 & 14.7499 \\ \hline
        32768  & 0.7387  & 0.9278  & 39.3491 & 10.7836 \\ \hline
        131072 & 0.2657  & 0.2101  & 39.2190 & 7.7920  \\ \hline
        \end{tabular}
    \end{table}

    \begin{table}[h!]
        \centering
        \caption{Example \ref{example2} -- $L^2$ norm errors for $\varepsilon=10^{-6}$ at $\Delta t = 0.001$ for ${\boldsymbol b}=(2,3)$.}
        \label{ex2-timeindependent-L2}
        \begin{tabular}{|l|l|l|l|l|}
        \hline
        $N$    & MM-FEM & MM-SUPG & FM-FEM & FM-SUPG \\ \hline
        512    & 0.2633 & 0.1289  & 1.5475 & 0.4138  \\ \hline
        2048   & 0.1127 & 0.0823  & 0.8708 & 0.2871  \\ \hline
        8192   & 0.0475 & 0.0096  & 0.5818 & 0.1573  \\ \hline
        32768  & 0.0019 & 0.0010  & 0.3720 & 0.0510  \\ \hline
        131072 & 0.0004 & 0.0005  & 0.1995 & 0.0180  \\ \hline
        \end{tabular}
    \end{table}

    Within this example, we also consider a time-dependent flow given by  ${\boldsymbol b} = (y-t, x-t)$. 
    This time-dependent velocity field introduces rotational and time-varying transport, testing robustness of the method under non-uniform flow. 
    The convergence rates in the $H^1$-seminorm and $L^2$ norm can be found in Figure \ref{ex2-timedependent}, 
    where numerical values are recorded in Tables \ref{ex2-timedependent-H1} and \ref{ex2-timedependent-L2}. Again, both MM-FEM and MM-SUPG perform much better than FM-FEM and FM-SUPG with larger $N$.

    \begin{figure}[thb]
        \centering
        \hbox{
        \begin{minipage}[t]{3.0in}
        \includegraphics[width=3.0in]{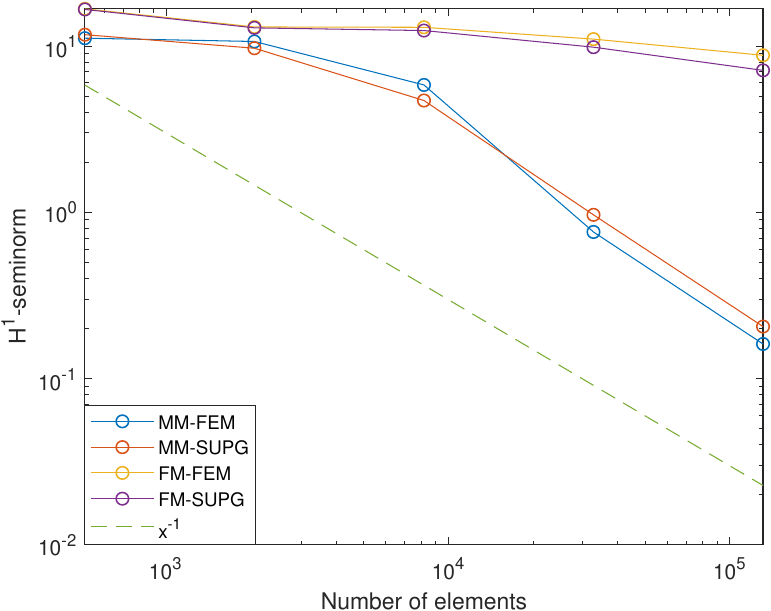}
        \centerline{(a): $H^1$-seminorm error.}
        \end{minipage}
        \hspace{10mm}
        \begin{minipage}[t]{3.0in}
        \includegraphics[width=3.0in]{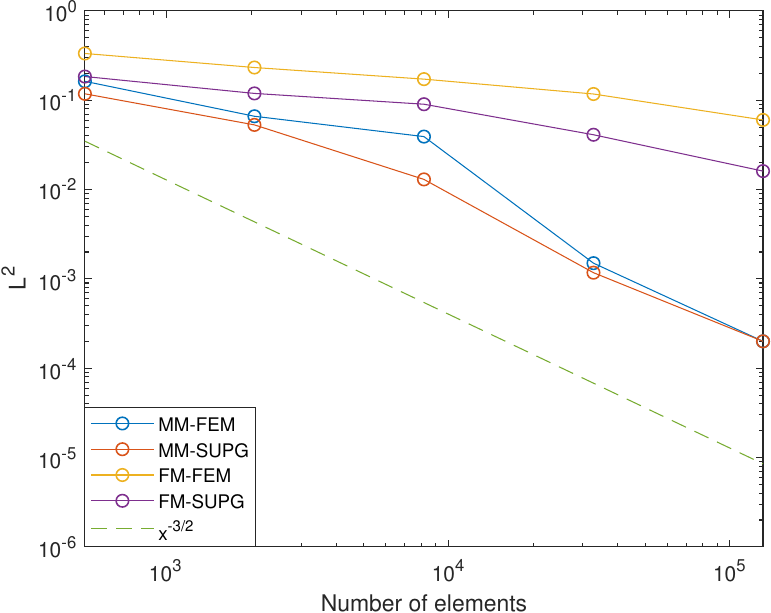}
        \centerline{(b): $L^2$ norm error.}
        \end{minipage}
        }
        \caption{Example \ref{example2} -- Errors corresponds to Tables \ref{ex2-timedependent-H1} and \ref{ex2-timedependent-L2}. $H^1$-seminorm and $L^2$ norm errors with $\varepsilon = 10^{-6}$ and $\Delta t = 0.001$ with flow field ${\boldsymbol b} = (y-t,x-t)$ at $T=0.5$.} 
        \label{ex2-timedependent}
    \end{figure}

    \begin{table}[thb]
        \centering 
        \caption{Example \ref{example2} -- $H^1$-seminorm errors for $\varepsilon = 10^{-6}$ at $\Delta t = 0.001$ for ${\boldsymbol b} = (y-t,x-t)$.}
        \label{ex2-timedependent-H1}
        \begin{tabular}{|l|l|l|l|l|}
        \hline
        $N$    & MM-FEM  & MM-SUPG & FM-FEM  & FM-SUPG \\ \hline
        512    & 11.1473 & 11.7017 & 16.7382 & 16.6045 \\ \hline
        2048   & 10.6529 & 9.7031  & 13.0265 & 12.8651 \\ \hline
        8192   & 5.8337  & 4.7036  & 12.9658 & 12.3947 \\ \hline
        32768  & 0.7614  & 0.9654  & 11.0075 & 9.8446  \\ \hline
        131072 & 0.1614  & 0.2054  & 8.7996  & 7.1405  \\ \hline
        \end{tabular}
    \end{table}

    \begin{table}[thb]
        \centering 
        \caption{Example \ref{example2} -- $L^2$ errors for $\varepsilon = 10^{-6}$ at $\Delta t = 0.001$ for ${\boldsymbol b} = (y-t,x-t)$.}
        \label{ex2-timedependent-L2}
        \begin{tabular}{|l|l|l|l|l|}
        \hline
        $N$    & MM-FEM   & MM-SUPG  & FM-FEM & FM-SUPG \\ \hline
        512    & 0.1629   & 0.1181   & 0.3340 & 0.1843  \\ \hline
        2048   & 0.0662   & 0.0530   & 0.2326 & 0.1193  \\ \hline
        8192   & 0.0392   & 0.0130   & 0.1724 & 0.0903  \\ \hline
        32768  & 0.0015   & 0.0011   & 0.1175 & 0.0411  \\ \hline
        131072 & 0.0002   & 0.0002   & 0.0602 & 0.0161  \\ \hline
        \end{tabular}
    \end{table}
    
    We see small oscillations along the top of the surface in the direction of the outflow in Figures \ref{ex2-sideview2}, even for large $N$. 
    However, for larger $N$, we see that DMP is not violated. Both methods converge to the same approximate solution as $N$ increases. 
    In Figure \ref{ex2-sideview} we see that the MM-SUPG method is a noticeable improvement compared to all methods. 
    In Figure \ref{ex2-sideview2}, only slight oscillations are observable on the backside of the surface (in the direction of the streamlines) for MM-FEM. This is more prominent for smaller $N$ in the MM-FEM case.
    The different results are due to the fact that Example \ref{example1} involves strong advective transport of sharp features, where SUPG stabilization plays a dominant role. 
    In contrast, Example \ref{example2} has a smoother structure where mesh adaptation alone already captures the solution well, reducing the relative impact of SUPG.
    Overall, MM-SUPG provides improved stability and reduced oscillations compared to fixed mesh methods, while offering comparable accuracy to MM-FEM.

    \begin{figure}[thb]
        \centering
        \hbox{
        \begin{minipage}[t]{3.0in}
        \includegraphics[width=3.0in]{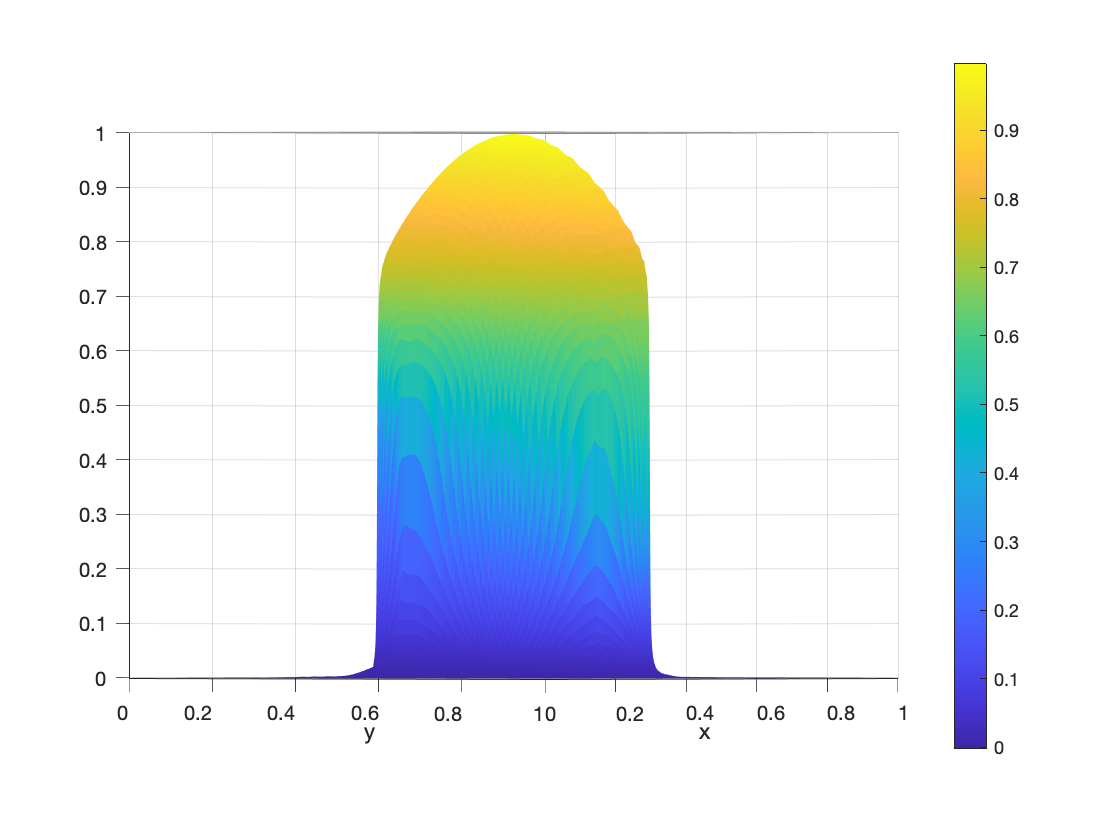}
        \centerline{(a): MM-FEM}
        \end{minipage}
        \hspace{10mm}
        \begin{minipage}[t]{3.0in}
        \includegraphics[width=3.0in]{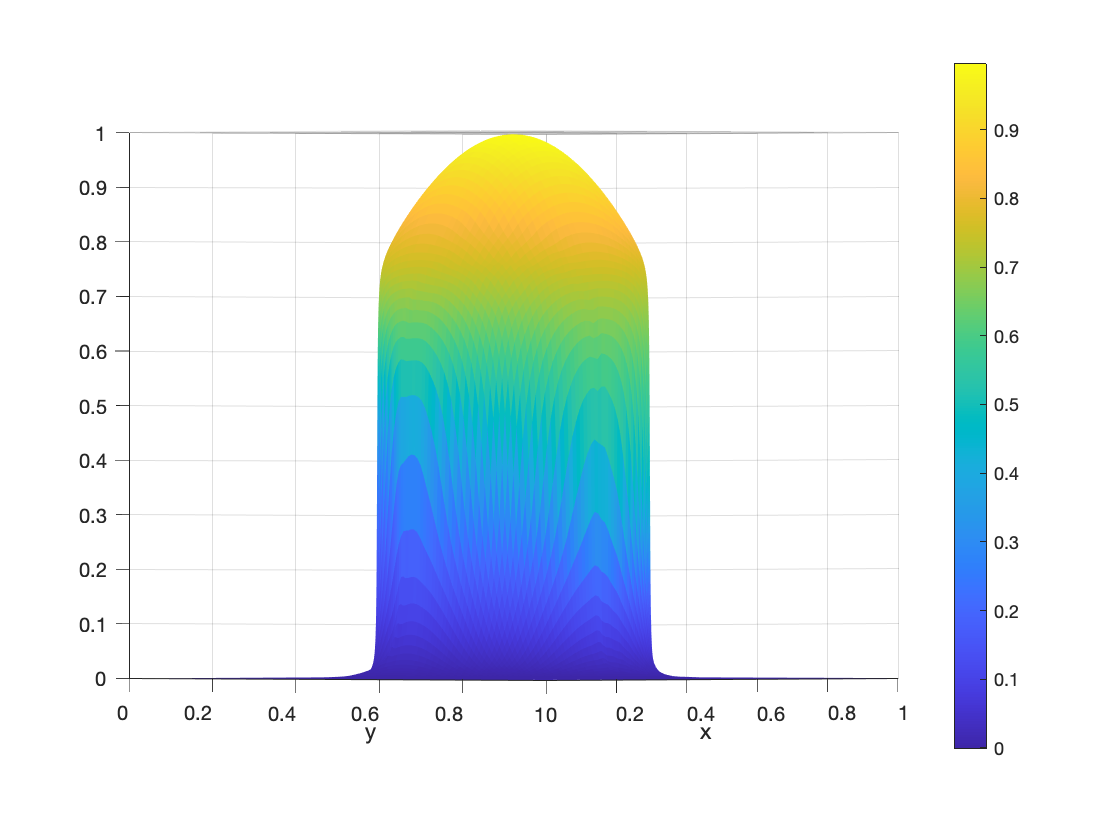}
        \centerline{(b): MM-SUPG}
        \end{minipage}
        }
        \vspace{0.5in}
        \hbox{
        \begin{minipage}[t]{3.0in}
        \includegraphics[width=3.0in]{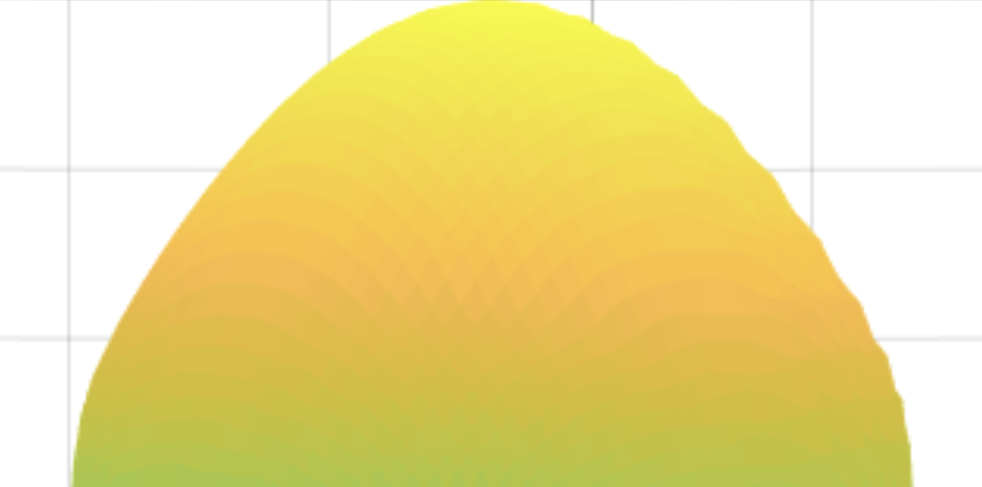}
        \centerline{(c): Enlarged view of (a)}
        \end{minipage}
        \hspace{10mm}
        \begin{minipage}[t]{3.0in}
        \includegraphics[width=3.0in]{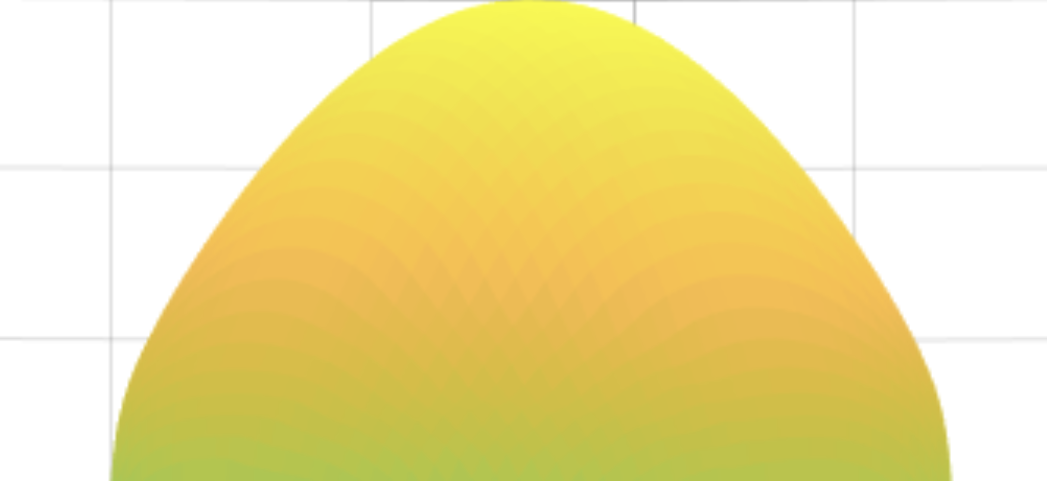}
        \centerline{(d): Enlarged view of (b)}
        \end{minipage}
        }
        \caption{Example \ref{example2} -- side view with $\Delta t = 0.001$ and ${\boldsymbol b} = (y-t,x-t)$, $\varepsilon = 10^{-6}$, and $N = 32768$ for MM-FEM and MM-SUPG at $T=0.5$.}
        \label{ex2-sideview2}
    \end{figure}

\end{exam}

\section{Discrete Maximum Principle for Anisotropic Diffusion}

In this section, we consider the anisotropic diffusion case where the diffusivity tensor $\mathbb{D}$ has different eigenvalues, 
which indicates diffusion is faster in one direction (the principal diffusion direction) than in the other direction. 

In the isotropic case, mesh adaptation primarily improves resolution of steep gradients. In contrast, for anisotropic diffusion, the orientation of the mesh relative to the diffusion tensor becomes critical. 
The discrete maximum principle is highly sensitive to mesh geometry and coefficient alignment. In the MM-SUPG framework, both the stiffness matrix and stabilization terms depend on the evolving mesh.
Misalignment may lead to violation of the discrete maximum principle (DMP), even in stabilized formulations. Therefore, the design of metric tensors must account for both diffusion anisotropy and convection direction.
The proposed tensor and metric intersection strategy aim to balance the competing effects of convection and anisotropic diffusion.

We say the solution $u$ to (\ref{modeleqn}) satisfies the weak maximum principle if $u\in C^{2,1}(\Omega_T)\cap C(\overline{\Omega}_T)$ and $f\leq 0$ in $\Omega_T$ \cite{evans1998partial}: 
\begin{align}
    \max_{({\boldsymbol x},t)\in \overline{\Omega}_T}u \leq \max\left\{0, \max_{({\boldsymbol x},t)\in \Gamma_T} u({\boldsymbol x},t)\right\},
\end{align}
where $\Gamma_T = \overline{\Omega}_T-\Omega_T$ is the parabolic boundary.

Let $S_i^K$ be the face opposite to the vertex ${\boldsymbol a}_i^K$, the $i$th vertex of element $K$. Introduce ${\boldsymbol q}$ vectors as \cite{li2013maximum}: 
\begin{align}
    \label{qvector}
    {\boldsymbol q}_i^K=\frac{{\boldsymbol n}_i^K}{h_i^K},
\end{align}
where ${\boldsymbol n}_i^K$ is the unit inward normal pointing to ${\boldsymbol a}_i^K$, and $h_i^K$ is the distance from vertex ${\boldsymbol a}_i^K$ to face $S_i^K$. 
Figure \ref{fig-qvectors} shows an illustration of ${\boldsymbol q}$ vectors.
\begin{figure}[thb]
    \centering
    \tikzset{every picture/.style={line width=0.75pt}} %set default line width to 0.75pt        
    \begin{tikzpicture}[x=0.75pt,y=0.75pt,yscale=-1,xscale=1]
        %Shape: Triangle [id:dp5673834222312759] 
        \draw  [fill={rgb, 255:red, 155; green, 155; blue, 155 }  ,fill opacity=0.65 ][line width=1.5]  (171.22,51.67) -- (160.97,104.94) -- (96.02,95.88) -- cycle ;
        %Shape: Polygon [id:dp34189967491494055] 
        \draw   (190.02,70.09) -- (140.31,134.56) -- (56.93,111.87) -- (55.1,33.37) -- (137.36,7.55) -- cycle ;
        %Straight Lines [id:da9223638972741786] 
        \draw    (55.1,33.37) -- (93.61,96.37) ;
        %Straight Lines [id:da25825226014827485] 
        \draw    (56.93,111.87) -- (93.61,96.37) ;
        %Straight Lines [id:da25587236123889623] 
        \draw    (137.36,7.55) -- (93.61,96.37) ;
        %Straight Lines [id:da525936205525142] 
        \draw    (93.57,96.53) -- (140.31,134.56) ;
        %Curve Lines [id:da5359694363119376] 
        \draw    (210.33,60) .. controls (249.57,30.94) and (289.29,30.51) .. (317.87,58.21) ;
        \draw [shift={(319.17,59.5)}, rotate = 225.32999999999998] [color={rgb, 255:red, 0; green, 0; blue, 0 }  ][line width=0.75]    (10.93,-3.29) .. controls (6.95,-1.4) and (3.31,-0.3) .. (0,0) .. controls (3.31,0.3) and (6.95,1.4) .. (10.93,3.29)   ;
        %Straight Lines [id:da5758999359075556] 
        \draw [line width=0.75]    (433.09,48.56) -- (405.52,83.41) ;
        \draw [shift={(404.28,84.98)}, rotate = 308.35] [fill={rgb, 255:red, 0; green, 0; blue, 0 }  ][line width=0.08]  [draw opacity=0] (12,-3) -- (0,0) -- (12,3) -- cycle    ;
        %Straight Lines [id:da40922161261584344] 
        \draw    (332.22,57.78) -- (373.3,74.99) ;
        \draw [shift={(375.15,75.77)}, rotate = 202.74] [fill={rgb, 255:red, 0; green, 0; blue, 0 }  ][line width=0.08]  [draw opacity=0] (12,-3) -- (0,0) -- (12,3) -- cycle    ;
        %Straight Lines [id:da6546177254305399] 
        \draw [line width=0.75]    (402,139.29) -- (402,83.23) ;
        \draw [shift={(402,81.23)}, rotate = 450] [fill={rgb, 255:red, 0; green, 0; blue, 0 }  ][line width=0.08]  [draw opacity=0] (12,-3) -- (0,0) -- (12,3) -- cycle    ;
        %Shape: Triangle [id:dp7196608138486933] 
        \draw   (369.98,29.44) -- (466.67,104.1) -- (337.34,104.1) -- cycle ;
        %Shape: Arc [id:dp8700393240299611] 
        \draw  [draw opacity=0] (344.56,87.1) .. controls (348.01,88.82) and (350.72,91.29) .. (352.39,94.49) .. controls (353.88,97.36) and (354.41,100.57) .. (354.08,103.93) -- (315.82,113.32) -- cycle ; \draw   (344.56,87.1) .. controls (348.01,88.82) and (350.72,91.29) .. (352.39,94.49) .. controls (353.88,97.36) and (354.41,100.57) .. (354.08,103.93) ;

        % Text Node
        \draw (135.5,79.07) node [anchor=north west][inner sep=0.75pt]    {$K$};
        % Text Node
        \draw (357.15,80.51) node [anchor=north west][inner sep=0.75pt]    {$\alpha _{ij}$};
        % Text Node
        \draw (378.87,107.4) node [anchor=north west][inner sep=0.75pt]    {$S_{i}^{K}$};
        % Text Node
        \draw (325.87,67.51) node [anchor=north west][inner sep=0.75pt]    {$S_{j}^{K}$};
        % Text Node
        \draw (385.33,63.73) node [anchor=north west][inner sep=0.75pt]    {$K$};
        % Text Node
        \draw (406.29,120.4) node [anchor=north west][inner sep=0.75pt]    {$\boldsymbol{q}_{i}^{K}$};
        % Text Node
        \draw (340.57,38.97) node [anchor=north west][inner sep=0.75pt]    {$\boldsymbol{q}_{j}^{K}$};
        % Text Node
        \draw (436.86,30.11) node [anchor=north west][inner sep=0.75pt]    {$\boldsymbol{q}_{r}^{K}$};
        % Text Node
        \draw (406.73,39.54) node [anchor=north west][inner sep=0.75pt]    {$S_{r}^{K}$};

    \end{tikzpicture}
    \caption{Illustration of faces and ${\boldsymbol q}$ vectors.}
    \label{fig-qvectors}
\end{figure}

In \cite{brandts2008discrete} it was shown that
\begin{align}
    \nabla\varphi_i\vert_K = {\boldsymbol q}_i^K.
\end{align}
Therefore, 
\begin{align}
    \label{energyintbasis}
    \int_{K}(\nabla\varphi_i)^2\dif {\boldsymbol x} =\|{\boldsymbol q}\|^2|K|=\frac{|K|}{(h_i^K)^2}.
\end{align}

Using the Divergence Theorem on an element $K$, we have
\begin{align}
    \int_K \varphi_i({\boldsymbol b}\cdot\nabla\varphi_i)\dif{\boldsymbol x}= \frac{1}{2}\int_{\partial K}\varphi_i^2({\boldsymbol b} \cdot {\boldsymbol n})\dif{\boldsymbol s} - \frac{1}{2}\int_{K}(\nabla\cdot{\boldsymbol b})\varphi_i^2\dif{\boldsymbol x}.
    \label{div-theorem}
\end{align}
Equation \eqref{div-theorem} does not determine the sign of the local convection term in general, since the boundary flux contribution on $\partial K$ need not vanish elementwise. 
Therefore, in the following analysis, we do not use \eqref{div-theorem} to infer positivity or negativity of $\int_K \varphi_i({\boldsymbol b}\cdot\nabla\varphi_i)\dif{\boldsymbol x}$. 
Instead, we use quantitative bounds as follows.
Using Cauchy--Schwarz inequality, we obtain
\begin{align}
\left|\int_K \phi_j (b \cdot \nabla \phi_i)\,dx\right|
&\le \|\phi_j\|_{L^2(K)} \|b \cdot \nabla \phi_i\|_{L^2(K)} \\
&\le \|\phi_j\|_{L^2(K)} \|b\|_{\infty,K} \|\nabla \phi_i\|_{L^2(K)}.
\end{align}
Since
\[
\int_K \phi_j^2\,dx = \frac{2|K|}{(d+1)(d+2)}
\quad \text{and} \quad
\|\nabla \phi_i\|_{L^2(K)}^2 = \frac{|K|}{(h_i^K)^2},
\]
it follows that
\begin{equation}
\left|\int_K \phi_j (b \cdot \nabla \phi_i)\,dx\right|
\le
\frac{\sqrt{2}\,|K|}{\sqrt{(d+1)(d+2)}\,h_i^K}\,\|b\|_{\infty,K}.
\end{equation}
A simpler estimate is obtained by noting that $\nabla \phi_i$ is constant on $K$, $\int_K \varphi_i \dif{\boldsymbol x} = \frac{K}{d+1}$, and $| \nabla \varphi_i| = \frac{1}{h_i^K}$:
\begin{align}
\left|\int_K \phi_j (b \cdot \nabla \phi_i)\,dx\right|
&\le \int_K |\phi_j|\,|b|\,|\nabla \phi_i|\,dx \\
\nonumber
&\le \|b\|_{\infty,K} |\nabla \phi_i| \int_K \phi_j\,dx = \frac{|K|}{(d+1)h_i^K}\,\|b\|_{\infty,K}.
\end{align}
Thus,
\begin{equation}
\label{convection-inequality}
\left|\int_K \phi_j (b \cdot \nabla \phi_i)\,dx\right| \le \frac{|K|}{(d+1)h_i^K}\,\|b\|_{\infty,K}.
\end{equation}

Similarly,
\begin{align}
\left|\int_K (b \cdot \nabla \phi_i)(b \cdot \nabla \phi_j)\,dx\right|
&\le \int_K |b|^2 |\nabla \phi_i| |\nabla \phi_j|\,dx \\
\nonumber
&\le \|b\|_{\infty,K}^2 |\nabla \phi_i| |\nabla \phi_j| |K| = \frac{|K|}{h_i^K h_j^K}\,\|b\|_{\infty,K}^2,
\end{align}
which gives
\begin{equation}
\label{diffusion-inequality}
\left|\int_K (b \cdot \nabla \phi_i)(b \cdot \nabla \phi_j)\,dx\right|
\le
\frac{|K|}{h_i^K h_j^K}\,\|b\|_{\infty,K}^2.
\end{equation}

For $\mathbb{D}_K$ in \eqref{galerkin}, following \cite{li2013maximum}, we define the diffusion-dependent map $\mathbb{D}_K^{-1/2}{\boldsymbol x}:K\to \tilde{K}$. 
Define the $\mathbb{D}_K$-norm by $\|{\boldsymbol x}\|_{\mathbb{D}_K}=\sqrt{{\boldsymbol x}^T\mathbb{D}_K{\boldsymbol x}}$. 
Then we have the following relations
\begin{align}
    \tilde{\boldsymbol a}_i^K = \mathbb{D}_K^{-1/2}{\boldsymbol a}_i^K,\ \ \tilde{S}_i^K=\mathbb{D}_K^{-1/2}S_i^K,\ \ |\tilde{K}|=\det(\mathbb{D}_K)^{1/2}|K|,\ \ \tilde{\boldsymbol q}_i^K=\mathbb{D}_K^{1/2}{\boldsymbol q}_i^K,\ \ \tilde{h}_i^K=\|\tilde{\boldsymbol q}_i^K\|^{-1}.
\end{align}
Denote by $\tilde{\alpha}_{ij}^K$ the dihedral angle between faces $\tilde{S}_i^K$ and $\tilde{S}_j^K$. Then 
\begin{align}
    \cos(\tilde{\alpha}_{ij}^K) = -\frac{(\tilde{\boldsymbol q}_i^K)^T\tilde{\boldsymbol q}_j^K}{\|\tilde{\boldsymbol q}_i^K\|\|\tilde{\boldsymbol q}_j^K\|}=-\frac{({\boldsymbol q}_i^K)^T\mathbb{D}_K{\boldsymbol q}_j^K}{\|{\boldsymbol q}_i^K\|_{\mathbb{D}_K}\|{\boldsymbol q}_j^K\|_{\mathbb{D}_K}},\ \ i\neq j.
\end{align}
Similarly, on the element $K$, we have the relation 
\begin{align}
    \cos(\alpha_{ij}) = -{\boldsymbol n}_i\cdot{\boldsymbol n}_j = -\frac{{\boldsymbol q}_i\cdot{\boldsymbol q}_j}{\|{\boldsymbol q}_i\|\|{\boldsymbol q}_j\|}.
\end{align}
It follows from (\ref{qvector}) that 
\begin{align}
    \label{eigenestimates}
    \frac{h_i^K}{\sqrt{\rho_{max}(\mathbb{D}_K})}\leq \tilde{h}_i^K\leq\frac{h_i^K}{\sqrt{\rho_{min}(\mathbb{D}_K})},
\end{align}
where $\rho_{min}(\mathbb{D}_K)$ and $\rho_{max}(\mathbb{D}_K)$ are the min and max eigenvalues of $\mathbb{D}_K$, respectively. 

The difficulty with the SUPG formulation \eqref{SUPGweakform} is the coupling of the diffusion tensor with the flow field ${\boldsymbol b}$ in \eqref{Amatrixint} as follows: 
\begin{align}
    (\mathbb{D}\nabla u_h,\nabla({\boldsymbol b}\cdot\nabla\varphi^h))_{0,K}.
\end{align}
For general ${\boldsymbol b} = (b_1(x_1,x_2),b_2(x_1,x_2))$, consider 
\begin{align}
    \label{dyadicform}
    \frac{\partial {\boldsymbol b}}{\partial x_i}\otimes\mathbb{D}_{(i)},\ \ i = 1,\dots, d,
\end{align}
where $\mathbb{D}_{(i)}$ is the vector $i$th row of $\mathbb{D}$. An essential assumption to the following analysis is that the tensor defined by 
\begin{align}
    \label{Kmatrix}
    \mathbb{K} = \sum_{i=1}^{d}\frac{\partial{\boldsymbol b}}{\partial x_i}\otimes\mathbb{D}_{(i)}, 
\end{align}
is symmetric and positive definite. This highlights a significant relationship between $\displaystyle\frac{\partial{\boldsymbol b}}{\partial x_i}$ and the $i$th row of $\mathbb{D}$ when using SUPG. 
For each element $K \in \mathcal{T}_h$, and denote the elementwise average of $\mathbb{K}$ as  
\begin{align}
    \label{K_K}
    \mathbb{K}_K = \frac{1}{|K|} \int_K \mathbb{K}(\boldsymbol{x}) d\boldsymbol{x},
\end{align}
according to \eqref{avg_tensor}. 

In the analysis below, we assume that $\mathbb{D}_K$ and $\mathbb{K}_K$ are symmetric and positive definite. 
This assumption imposes a structural alignment between the velocity field and the diffusion tensor. While restrictive, it enables tractable analysis of DMP preservation. 
This assumption is not expected to hold in general flows but serves as a sufficient condition.

The added diffusion integral term in (\ref{SUPGweakform}) can be written in general form as
\begin{align}
    \label{K_K_int}
    \tau_K\int_{K}(\mathbb{D}\nabla u)\cdot \nabla ({\boldsymbol b}\cdot\nabla v)\dif {\boldsymbol x} 
    &= \tau_K\sum_{i=1}^d\int_{K}(\nabla v)^T\left(\frac{\partial {\boldsymbol b}}{\partial x_i}\otimes\mathbb{D}_{(i)}\right)\nabla u \dif {\boldsymbol x}
    = \tau_K (\nabla v)^T \mathbb{K}_K \nabla u.
\end{align}

The analysis of the numerical scheme relies on the following fundamental lemma \cite{stoyan1986maximum}. 
\begin{lemma}
If $\tilde{\boldsymbol A}$ in \eqref{ODEsystem} is an $M$-matrix and has nonnegative row sums, then the numerical scheme satisfies DMP. 
\end{lemma}
The following lemmas were used in \cite{li2010anisotropic} to analyze DMP with ${\boldsymbol b} = {\boldsymbol 0}$. 
\begin{lemma}
\label{lemma-mmatrix}
If the mesh satisfies the anisotropic nonobtuse angle condition (ANAC)
\begin{align}
    \label{nonobtuse}
    0<\tilde{\alpha}_{ij}^K\leq \frac{\pi}{2},\ \ \forall i,j=1,\dots, d+1, i\neq j, \forall K\in\mathcal{T}_h,
\end{align}
then for any $\mathbb{F}_K$ of SPD-type, the matrix defined by 
\begin{align}
    \sum_{K\in\mathcal{T}_h}|K|(\nabla\varphi_i)^T\mathbb{F}_K\nabla\varphi_j,\ \  i =,1\dots N_{vi}, j=1,\dots, N_v
\end{align}
is an $M$-matrix and has nonnegative row sums. Moreover, for any element $K\in\mathcal{T}_h$, and $i,j = 1,\dots, d+1$, 
\begin{align}
\label{formulaforquadraticform}
(\nabla \varphi_i)^T\mathbb{F}_K\nabla\varphi_j = \begin{cases}
\displaystyle -\frac{\cos(\tilde{\alpha}_{ij}^K)}{\tilde{h}_i^K\tilde{h}_j^K}\textrm{\ \ for \ \ } i\neq j\\
\displaystyle \frac{1}{(\tilde{h}_i^K)^2}\textrm{\ \ for \ \ i = j}
\end{cases}.
\end{align}
When $d=2$, the case $i\neq j$ simplifies to 
\begin{align}
|K|(\nabla \varphi_i)^T\mathbb{F}_K\nabla\varphi_j = -\frac{\sqrt{\det(\mathbb{F}_K})}{2}\cot(\tilde{\alpha}_{ij}^K). 
\end{align}
\end{lemma}

We establish DMP in a series of steps. We first investigate spatial estimates, which are critical to understanding the required 
structure of the matrices ${\boldsymbol A}$ and ${\boldsymbol M}$. Temporal estimates are then derived in order to apply Lemma 4.1.

\subsection{General flow forms}
In this section we investigate the structure of $\mathbb{K}$ specified by ${\boldsymbol b}$. For $d = 2$, we have 
\begin{align}
    \label{tensorproduct1}
    \frac{\partial {\boldsymbol b}}{\partial x_1} \otimes \mathbb{D}_{(1)} =  \begin{pmatrix}
    \displaystyle\frac{\partial b_1}{\partial x_1}D_{11} & \displaystyle\frac{\partial b_1}{\partial x_1} D_{12}\\
    \displaystyle\frac{\partial b_2}{\partial x_1}D_{11} & \displaystyle\frac{\partial b_2}{\partial x_1}D_{12}
    \end{pmatrix},
\end{align}
and 
\begin{align}
    \label{tensorproduct2}
    \frac{\partial {\boldsymbol b}}{\partial x_2} \otimes \mathbb{D}_{(2)} =  \begin{pmatrix}
    \displaystyle\frac{\partial b_1}{\partial x_2}D_{21} & \displaystyle\frac{\partial b_1}{\partial x_2} D_{22}\\
    \displaystyle\frac{\partial b_2}{\partial x_2}D_{21} & \displaystyle\frac{\partial b_2}{\partial x_2}D_{22}
    \end{pmatrix}.
\end{align}
There are two ways in which we consider $\mathbb{K}$. If (\ref{tensorproduct1}) and (\ref{tensorproduct2}) are of SPD-type, then $\mathbb{K}$ will be as well. If not, then we need to analyze $\mathbb{K}$ in its general form given by 
\begin{align}
    \mathbb{K} = \begin{pmatrix}
    \frac{\partial b_1}{\partial x_1}D_{11}+\frac{\partial b_1}{\partial x_2}D_{21} & \frac{\partial b_1}{\partial x_1}D_{12}+\frac{\partial b_1}{\partial x_2}D_{22}\\
    \frac{\partial b_2}{\partial x_1}D_{11}+\frac{\partial b_2}{\partial x_2}D_{21} & \frac{\partial b_2}{\partial x_1}D_{12}+\frac{\partial b_2}{\partial x_2}D_{22}
    \end{pmatrix}.
\end{align}
Hence, using the symmetric condition, we would require
\begin{align}
    \frac{\partial}{\partial x_1}(b_2D_{11}-b_1D_{12}) = \frac{\partial }{\partial x_2}(b_1D_{22} - b_2D_{21}). 
\end{align}
Assume first that these matrices are of SPD-type so that general form of $\mathbb{K}$ is of SPD-type.

\begin{theorem}
Suppose $\mathbb{D}$, equations (\ref{tensorproduct1}), and (\ref{tensorproduct2}) are of SPD-type with constant real-valued entries and $D_{12} = D_{21}\neq 0$. Then there are functions $f_1(x_2,t)$, $g_1(x_2,t)$, $f_2(x_1,t)$, and $g_2(x_1,t)$ such that 
\begin{align}
    b_1(x_1,x_2,t) = \frac{D_{11}D_{22}f_1(x_1,t)-D_{12}D_{21}f_2(x_2,t)}{\det(\mathbb{D})},
\end{align}
and 
\begin{align}
    b_2(x_1,x_2,t)& =\frac{D_{12}D_{22}(f_1(x_1,t) - f_2(x_2,t))}{\det(\mathbb{D})},
\end{align}
for $D_{22} \neq 0$. Likewise, 
\begin{align}
    b_1(x_1,x_2,t) = \frac{D_{11}D_{21}(g_2(x_2,t) - g_1(x_1,t))}{\det(\mathbb{D})},
\end{align}
and 
\begin{align}
    b_2(x_1,x_2,t) = \frac{D_{11}D_{22}g_2(x_2,t)-D_{12}D_{21}g_1(x_1,t)}{\det(\mathbb{D})},
\end{align}
whenever $D_{11}\neq 0$. 
\end{theorem}
\begin{proof}
Assume $D_{12}=D_{21}\neq 0$ and $D_{22}\neq 0$. The other case is very similar. Assuming symmetry as mentioned, we have 
\begin{align}
    \frac{\partial b_1}{\partial x_1} &= \left(\frac{D_{11}}{D_{12}}\right)\frac{\partial b_2}{\partial x_1}, \\
    b_1(x_1,x_2,t) &= \left(\frac{D_{11}}{D_{12}}\right)b_2(x_1,x_2,t) + f_2(x_2,t),
\end{align}
for some $f_2(x_2,t)$. Additionally, 
\begin{align}
    \frac{\partial b_1}{\partial x_2} &= \left(\frac{D_{21}}{D_{22}}\right)\frac{\partial b_2}{\partial x_2}, \\
   b_1(x_1,x_2,t) &= \left(\frac{D_{21}}{D_{22}}\right)b_2(x_1,x_2,t) + f_1(x_1,t),
\end{align}
for some function $f_1(x_1,t)$. Equating the two equations for $b_1$ gives 
\begin{align}
    \left(\frac{D_{11}}{D_{12}}\right)b_2(x_1,x_2,t) + f_2(x_2,t) &= \left(\frac{D_{21}}{D_{22}}\right)b_2(x_1,x_2,t) + f_1(x_1,t), \\
    \left(\frac{D_{11}}{D_{12}} -\frac{D_{21}}{D_{22}}\right)b_2(x_1,x_2,t) &= f_1(x_2,t) - f_2(x_1,t), \\
    b_2(x_1,x_2,t)& =\frac{D_{12}D_{22}(f_1(x_1,t) - f_2(x_2,t))}{\det(\mathbb{D})}.
\end{align}
Plugging this into the equation for $b_1$, we have 
\begin{align}
    \nonumber
     b_1(x_1,x_2,t) &= \left(\frac{D_{11}}{D_{12}}\right)\frac{D_{12}D_{22}(f_1(x_1,t) - f_2(x_2,t))}{\det(\mathbb{D})} + f_2(x_2,t)\\
     \nonumber
     & = \frac{D_{11}D_{22}(f_1(x_1,t)-f_2(x_2,t))+f_2(x_2,t)(D_{11}D_{22}-D_{12}D_{21})}{\det(\mathbb{D})}\\
     & = \frac{D_{11}D_{22}f_1(x_1,t)-D_{12}D_{21}f_2(x_2,t)}{\det(\mathbb{D})}.
\end{align}
Note that since $\mathbb{D}$ is of SPD-type, $\det(\mathbb{D})\neq 0$. 
\end{proof}

\begin{rem}
\label{rem:generalflow}
If $D_{12}=D_{21}$ is not zero, the general form of ${\boldsymbol b}$ is given by 
\begin{align}
    \frac{1}{\det(\mathbb{D})}\left(D_{22}f_1(x_1,t)\mathbb{D}_{(1)}-D_{12}f_2(x_2,t)\mathbb{D}_{(2)}\right),
\end{align}
for some functions $f_1$ and $f_2$ whenever $D_{22}\neq 0$. Similarly, if $D_{11}\neq 0$, then the general form of ${\boldsymbol b}$ is given by 
\begin{align}
    \frac{1}{\det(\mathbb{D})}\left(D_{11}g_2(x_2,t)\mathbb{D}_{(2)}-D_{21}g_1(x_1,t)\mathbb{D}_{(1)}\right).
\end{align}
If $D_{11}$ and $D_{22}$ are both not zero, then $g_1 = -\displaystyle\frac{D_{22}}{D_{21}}f_1$ and $g_2 = -\displaystyle\frac{D_{12}}{D_{11}}f_2$.
\end{rem}

\begin{corollary}
If $\mathbb{D} = \textrm{diag}(D_{11},D_{22})$, with $D_{11}$ and $D_{22}\neq 0$, then there are functions $f$ and $g$ such that 
\begin{align}
    \label{bform2}
    {\boldsymbol b}(x_1,x_2,t) = \begin{pmatrix} \displaystyle\frac{f(x_1,t)}{D_{11}}\\
    \displaystyle\frac{g(x_2,t)}{D_{22}}
    \end{pmatrix},
\end{align}
and 
\begin{align}
    \mathbb{K} = \begin{pmatrix}
    \displaystyle \frac{\partial f}{\partial x_1} & 0 \\
    0 & \displaystyle\frac{\partial g}{\partial x_2}
    \end{pmatrix}.
\end{align}
\end{corollary}
\begin{proof}
The form of ${\boldsymbol b}$ as shown in (\ref{bform2}) follows from a similar argument as shown in the proof of the previous theorem. 
\end{proof}

Based on the Remark \ref{rem:generalflow}, we require 
\begin{align}
    \frac{D_{22}}{\det(\mathbb{D})}\left(D_{11}\frac{\partial f_1}{\partial x_1}-D_{12}\frac{\partial f_2}{\partial x_2}\right)
\end{align}
to share the same sign as the divergence of ${\boldsymbol b}$. 

\subsection{Mass and stiffness matrix structure: estimates on $\tau_K$}
We first establish non-negative diagonal entries for the mass and stiffness matrices. 

\begin{lemma}
\label{lemma-diagonal}
Suppose $\mathbb{K}_K$ and $\mathbb{D}_K$ are symmetric and positive definite and the meshes associated with $\mathbb{D}_K$ and $\mathbb{K}_K$ satisfy ANAC. 
Then the diagonal entries of the mass and stiffness matrices satisfy $M_{ii} > 0$ and $A_{ii} > 0$ provided
\begin{equation}
\tau_K \|b\|_{\infty,K} \le \frac{h_i^K}{d+2},
\label{miicondition}
\end{equation}
and
\begin{equation}
\frac{\rho_{\min}( \mathbb{D}_K) + \tau_K \rho_{\min}(\mathbb{K}_K)}{(h_i^K)^2} > \frac{1}{2\varepsilon\tau_K (d+1)(d+2)}.
\label{aiicondition}
\end{equation}
\end{lemma}
\begin{proof}
We consider a patch of elements $\omega_i$ that share a vertex ${\boldsymbol a}_i$.
Using \eqref{convection-inequality} and \eqref{basisformula}, from \eqref{MmatrixFormula}, we obtain
\[
M_{ii}
\ge
\sum_{K \in \omega_i}
\left(
\frac{2|K|}{(d+1)(d+2)}
-
\tau_K \frac{|K|}{(d+1)h_i^K}\|b\|_{\infty,K}
\right),
\]
which is positive under \eqref{miicondition}.

Applying Young's inequality
\[
\phi_i (b \cdot \nabla \phi_i) + \tau_K (b \cdot \nabla \phi_i)^2
\ge -\frac{1}{4\tau_K}\phi_i^2,
\]
for $A_{ii}$ in \eqref{Amatrixint},
we have,
\begin{align}
    \label{aiiformula}
    \nonumber
    A_{ii}&= \sum_{K\in\omega_i}\int_K \left( \varphi_i{\boldsymbol b}\cdot \nabla\varphi_i+\tau_K({\boldsymbol b\cdot \nabla\varphi_i})^2+\tau_K\left(\varepsilon \mathbb{D}_K \nabla\varphi_i\right)\cdot\nabla({\boldsymbol b\cdot\nabla \varphi_i)} \right) \dif {\boldsymbol x}+|K|(\nabla\varphi_i)^T \varepsilon \mathbb{D}_K\nabla\varphi_i \\
    \nonumber
    &= \sum_{K\in\omega_i}\int_K \left( \varphi_i{\boldsymbol b}\cdot \nabla\varphi_i+\tau_K({\boldsymbol b\cdot \nabla\varphi_i})^2  \right) \dif {\boldsymbol x} + \tau_K |K|(\nabla\varphi_i)^T \varepsilon \mathbb{K}_K\nabla\varphi_i + |K|(\nabla\varphi_i)^T \varepsilon \mathbb{D}_K\nabla\varphi_i \\
    \nonumber
    & \ge \sum_{K\in\omega_i} \int_K -\frac{1}{4\tau_K}\phi_i^2 \dif {\boldsymbol x}  + |K| (\nabla\varphi_i)^T \varepsilon (\tau_K \mathbb{K}_K + \mathbb{D}_K) \nabla\varphi_i \\
    \nonumber
    & = \sum_{K\in\omega_i} -\frac{|K|}{2\tau_K (d+1)(d+2)} + |K| (\nabla\varphi_i)^T \varepsilon (\tau_K \mathbb{K}_K + \mathbb{D}_K) \nabla\varphi_i \\
    & \ge \sum_{K\in\omega_i} |K| \left( \varepsilon \frac{\rho_{min}(\mathbb{D}_K)+\tau_K \rho_{min}(\mathbb{K}_K)}{(h_i^K)^2} - \frac{1}{2\tau_K (d+1)(d+2)}\right), 
\end{align}
which yields positivity under \eqref{aiicondition}.
\end{proof}

It is important to keep in mind that this estimate is only for those elements associated with $A_{ii}$ on the interior of $\Omega$. This estimate highlights the significant relationship between $\tau_K$ and the convection term. 

\begin{lemma}
\label{lemma-aij}
Suppose $\mathbb{K}_K$ and $\mathbb{D}_K$ are symmetric and positive definite and the meshes associated with $\mathbb{D}_K$ and $\mathbb{K}_K$ satisfy ANAC.
For $i \ne j$, define
\begin{equation}
\label{delta_ij}
\delta_{ij}^K =
\rho_{\min}(\varepsilon \mathbb{D}_K)\cos(\widetilde{\alpha}_{ij}^{\mathbb{D}}) + \tau_K \rho_{\min}(\mathbb{K}_K)\cos(\check{\alpha}_{ij}^\mathbb{K}) 
- \frac{h_i^K}{d+1}\|b\|_{\infty,K} - \tau_K \|b\|_{\infty,K}^2,
\end{equation}
where $\widetilde{\alpha}_{ij}^{\mathbb{D}}$ and $\widetilde{\alpha}_{ij}^\mathbb{K}$ are the angles in the meshes associated with $\mathbb{D}_K$ and $\mathbb{K}_K$, respectively. 
If 
\begin{equation}
    \label{delta_ij_cond}
    \delta_{ij}^K \ge 0 \quad \text{ for all } K \in \mathcal{T}_h, \, i \neq j,
\end{equation}
then $A_{ij} \le 0$.
\end{lemma}

\begin{proof}
For $i \ne j$, from \eqref{Amatrixint} we have 
\begin{align}
    \nonumber
    A_{ij} = \sum_{K\in\omega_i\cap\omega_j} \left( \int_K \varphi_i({\boldsymbol b}\cdot \nabla\varphi_j)\dif {\boldsymbol x} +\tau_K \int_K ({\boldsymbol b\cdot \nabla\varphi_j})({\boldsymbol b\cdot \nabla\varphi_i}) \dif {\boldsymbol x} 
    + |K|(\nabla\varphi_i)^T (\varepsilon \mathbb{D}_K + \varepsilon \tau_K \mathbb{K}_K)\nabla\varphi_j \right).
\end{align}

Applying Lemma \ref{lemma-mmatrix} to $\varepsilon \mathbb{D}_K$ and $\varepsilon \mathbb{K}_K$, we have
\begin{equation}
    |K| \nabla\varphi_i)^T (\varepsilon \mathbb{D}_K)\nabla\varphi_j 
    \le -\frac{|K|}{h_i^K h_j^K} \rho_{min}(\varepsilon \mathbb{D}_K)\cos(\widetilde{\alpha}_{ij}^{\mathbb{D}}), 
\end{equation}
and
\begin{equation}
    |K| \nabla\varphi_i)^T (\varepsilon \mathbb{K}_K)\nabla\varphi_j 
    \le -\frac{|K|}{h_i^K h_j^K} \rho_{min}(\varepsilon \mathbb{K}_K)\cos(\widetilde{\alpha}_{ij}^\mathbb{K}). 
\end{equation}
Using inequalities \eqref{convection-inequality} and \eqref{diffusion-inequality}, we have
\begin{align}
    \nonumber
    A_{ij} &\le \sum_{K\in\omega_i\cap\omega_j} \frac{|K|}{h_i^K h_j^K} \left( \frac{h_i^K}{d+1} \|b\|_{\infty,K} + \tau_K \|b\|^2_{\infty,K} 
    - \rho_{min}(\varepsilon \mathbb{D}_K)\cos(\widetilde{\alpha}_{ij}^{\mathbb{D}}) - \tau_K \rho_{min}(\varepsilon \mathbb{K}_K)\cos(\widetilde{\alpha}_{ij}^{\mathbb{K}})  \right) \\
    \label{aij-condition}
    &= - \sum_{K\in\omega_i\cap\omega_j} \frac{|K|}{h_i^K h_j^K} \delta_{ij}^K.
\end{align}
Therefore, $A_{ij} \le 0$ if \eqref{delta_ij_cond} holds.
\end{proof}
In real computations, there is only one mesh at each physical time, thus $\widetilde{\alpha}_{ij}^\mathbb{D}$ and $\widetilde{\alpha}_{ij}^\mathbb{K}$ will be referred to the same angle, 
which will be discussed in later subsection \ref{subsec-metric}.

\begin{lemma}
\label{lemma-mij}
For $i \ne j$, if
\begin{align}
    \label{cond-mij}
    \tau_K \|b\|_{\infty,K} \le \frac{h_i^K}{d+2} \quad \text{ for all } K \in \mathcal{T}_h, 
\end{align}
then $M_{ij}\geq 0$.
\end{lemma}
\begin{proof}
For $i \ne j$, from \eqref{MmatrixFormula} and \eqref{basisformula}
\begin{equation}
M_{ij} = \sum_{K\in\omega_i\cap\omega_j} \left( \frac{|K|}{(d+1)(d+2)} + \tau_K \int_K \varphi_j({\boldsymbol b}\cdot \nabla\varphi_i)\dif {\boldsymbol x} \right).
\end{equation}
Applying \eqref{convection-inequality}, we have
\begin{equation}
M_{ij} \ge \sum_{K\in\omega_i\cap\omega_j} \left( \frac{|K|}{(d+1)(d+2)} - \tau_K \frac{|K|}{(d+1)h_i^K}\,\|b\|_{\infty,K} \right),
\end{equation}
which is nonnegative under condition \eqref{cond-mij}.
\end{proof}

Combining the previous two lemmas, we have the following result immediately. 
\begin{lemma}
Assume the hypotheses of Lemmas \ref{lemma-aij} and \ref{lemma-mij}. Then
\begin{align}
    \label{nonnegativecondition}
    M_{ij} - \Delta t(1-\theta)A_{ij}\geq 0, \quad i \ne j.
\end{align}
\end{lemma}

\subsection{Temporal estimates}
In this section we derive upper and lower bounds on the fixed time step $\Delta t$ so that DMP is satisfied. 

Define 
\begin{align}
    \eta_{ij}^K := \frac{|K|}{h_i^k h_j^K} \delta_{ij}^K, \quad \mu_{ij}^K := \frac{|K|}{(d+1)(d+2)} + \tau_K \frac{|K|}{(d+1)h_i^K} \|b\|_{\infty,K}.
\end{align}
\begin{lemma}
Assume that for each element $K$ in $\mathcal{T}_h$:
\begin{itemize}
    \item [(i)] The elementwise tensors $\mathbb{D}_K$ and $\mathbb{K}_K$ are symmetric and positive definite.
    \item [(ii)] The meshes associated with $\mathbb{D}_K$ and $\mathbb{K}_K$ satisfy the anisotropic nonobtuse angle condition (ANAC).
    \item [(iii)] The stabilization parameter $\tau_K$ satisfies 
    \begin{equation}
        \label{cond-stab-1}
        \tau_K \|b\|_{\infty,K} \le C_1 h_K,
    \end{equation}
    and 
    \begin{equation}
        \label{cond-stab-2}
        \varepsilon (\rho_{min}(\mathbb{D}_K)+\tau_K \rho_{min}(\mathbb{K}_K)) > C_2 \frac{h^2_K}{\tau_K},
    \end{equation}
    for positive constants $C_1$, and $C_2$ independent of the mesh. 
    \item [(iv)] For all $i \ne j$, $\delta_{ij}^K$ defined in \eqref{delta_ij} is positive. 
\end{itemize}
Then, if 
\begin{equation}
\label{temporallowerbound}
\Delta t \ge \frac{1}{\theta} \max_{i \ne j} \frac{\sum_{K\in\omega_i\cap\omega_j} \mu_{ij}^K}{\sum_{K\in\omega_i\cap\omega_j}\eta_{ij}^K},
\end{equation}
the matrix $\tilde{\boldsymbol A} = {\boldsymbol M} + \theta \Delta t {\boldsymbol A}$ is an $M$-matrix. 
\end{lemma}
\begin{proof}
We first show that $\tilde{\boldsymbol A}$ is a $Z$-matrix. 
For $i \ne j$, from \eqref{MmatrixFormula} and \eqref{basisformula}
\begin{equation}
M_{ij} = \sum_{K\in\omega_i\cap\omega_j} \left( \frac{|K|}{(d+1)(d+2)} + \tau_K \int_K \varphi_j({\boldsymbol b}\cdot \nabla\varphi_i)\dif {\boldsymbol x} \right).
\end{equation}
Applying \eqref{convection-inequality}, we have
\begin{equation}
M_{ij} \le \sum_{K\in\omega_i\cap\omega_j} \left( \frac{|K|}{(d+1)(d+2)} + \tau_K \frac{|K|}{(d+1)h_i^K}\,\|b\|_{\infty,K} \right) 
= \sum_{K\in\omega_i\cap\omega_j} \mu_{ij}^K,
\end{equation}
and from \eqref{aij-condition}, we obtain
\begin{equation}
A_{ij} \le - \sum_{K\in\omega_i\cap\omega_j} \eta_{ij}^K.
\end{equation}
Thus,
\begin{equation}
\tilde{\boldsymbol A}_{ij} = M_{ij} + \theta \Delta t A_{ij} \le \sum_{K\in\omega_i\cap\omega_j} \mu_{ij}^K - \theta \Delta t \sum_{K\in\omega_i\cap\omega_j} \eta_{ij}^K \le 0
\end{equation}
provided \eqref{temporallowerbound} holds. Hence all off-diagonal entries of $\tilde{\boldsymbol A}$ are nonpositive. 
Combined with the positivity of the diagonal entries of ${\boldsymbol M}$ and ${\boldsymbol A}$, $\tilde{\boldsymbol A}$ is a $Z$-matrix. 

Next, we show that $\tilde{\boldsymbol{A}}$ is an $M$-matrix. 
Reorder the system into interior and boundary nodes, by (\ref{Amatrix}), we can write 
\begin{align*}
    \tilde{\boldsymbol A} = \begin{bmatrix}
    \boldsymbol{B} & \boldsymbol{C}\\
    0 & \boldsymbol{I}
    \end{bmatrix}, \quad \boldsymbol{B} := \boldsymbol{M}_{11} + \theta\Delta t \boldsymbol{A}_{11}. 
\end{align*}
Since $\boldsymbol{C} \le 0$, it suffices to show that $\boldsymbol{B}$ is an $M$-matrix. 

Firstly, the block $\boldsymbol{B}$ is also a $Z$-matrix because $\tilde{\boldsymbol A}$ is a $Z$-matrix.  
Next, we show $\boldsymbol{B}$ is positive-definite. 

Let ${\boldsymbol v} = (v_1, v_2, \dots, v_{N_{vi}})$ be a nonzero vector in $\mathbb{R}^{N_{vi}}$ and write $v_h = \sum_{j=1}^{N_{vi}}v_j\varphi_j \in U_0^h$. 
Then, 
\begin{equation}
    \boldsymbol{v}^T \boldsymbol{B} \boldsymbol{v} = \boldsymbol{v}^T \boldsymbol{M}_{11} \boldsymbol{v} + \theta\Delta t \boldsymbol{v}^T \boldsymbol{A}_{11} \boldsymbol{v},
\end{equation}
and the interior mass bloack satisfies
\begin{equation}
    \boldsymbol{v}^T \boldsymbol{M}_{11} \boldsymbol{v} = \sum_{i,j=1}^{N_{vi}} v_i M_{ij} v_j. 
\end{equation}
Using the definition of $M_{ij}$ in \eqref{MmatrixFormula}, this becomes
\begin{equation}
    \boldsymbol{v}^T \boldsymbol{M}_{11} \boldsymbol{v} = \sum_{K \in \mathcal{T}_h} \left( \int_K v_h^2 \dif \boldsymbol{x} + \tau_K \int_K v_h (\boldsymbol{b} \cdot \nabla v_h) \dif \boldsymbol{x} \right). 
\end{equation}
Applying Cauchy-Schwarz and inverse inequality, we have
\begin{equation}
    \left| \int_K v_h (\boldsymbol{b} \cdot \nabla v_h) \dif \boldsymbol{x} \right| \le C_{inv} h_K^{-1} \|b\|_{\infty,K} \|v_h\|^2_{L^2(K)}.
\end{equation}
Therefore,
\begin{align}
    \nonumber
    \boldsymbol{v}^T \boldsymbol{M}_{11} \boldsymbol{v} &\ge  \sum_{K \in \mathcal{T}_h} \left( \|v_h\|^2_{L^2(K)} - \tau_K C_{inv} h_K^{-1} \|b\|_{\infty,K} \|v_h\|^2_{L^2(K)} \right) \\
    \nonumber
    & = \sum_K \left(1 - C \tau_K h_K^{-1} \|b\|_{\infty,K} \right)\|v_h\|^2_{L^2(K)} \\
    \nonumber
    & > 0,
\end{align}
provided $\tau_K$ satisfies the condition \eqref{cond-stab-1}.

Similarly, from \eqref{Amatrixint} and \eqref{K_K_int}, we get,
\begin{align}
    \label{v_A_v}
    \boldsymbol{v}^T \boldsymbol{A}_{11} \boldsymbol{v} &= \sum_{K \in \mathcal{T}_h} \left( \int_K v_h (\boldsymbol{b} \cdot \nabla v_h) \dif \boldsymbol{x} 
    + \tau_K \int_K (\boldsymbol{b} \cdot \nabla v_h)^2 \dif \boldsymbol{x} \right) \\
    \nonumber
    &+ \sum_{K \in \mathcal{T}_h} \varepsilon |K| (\nabla v_h)^T (\mathbb{D}_K + \tau_K \mathbb{K}_K) \nabla v_h. 
\end{align}
Applying Young's inequality to the scaler quantities $v_h$ and $\boldsymbol{b} \cdot \nabla v_h$, we have
\begin{equation}
    v_h (b \cdot \nabla \phi_i) + \tau_K (b \cdot \nabla \phi_i)^2 \ge -\frac{1}{4\tau_K} v_h^2.
\end{equation}
Integrate over $K$ and substitue into \eqref{v_A_v}, we have
\begin{align}
    \nonumber
    \boldsymbol{v}^T \boldsymbol{A}_{11} \boldsymbol{v} &\ge \sum_{K \in \mathcal{T}_h} \left( - \frac{1}{4\tau_K} \|v_h\|^2_{L^2(K)}
    + \varepsilon |K| (\nabla v_h)^T (\mathbb{D}_K + \tau_K \mathbb{K}_K) \nabla v_h \right). 
\end{align}

Using SPD properties of $\mathbb{D}_K$ and $\mathbb{K}_K$, we obtain
\begin{equation}
    (\nabla v_h)^T (\mathbb{D}_K + \tau_K \mathbb{K}_K) \nabla v_h \ge (\rho_{min}(\mathbb{D}_K) + \tau_K \rho_{min}(\mathbb{K}_K)) |\nabla v_h|^2.
\end{equation}
Since $v_h$ is linear on each element $K$, $\nabla v_h$ is constant on $K$, so
\begin{equation}
    |K| |\nabla v_h|^2 = \|\nabla v_h\|^2_{L^2(K)}.
\end{equation}
Therefore, 
\begin{align}
    \nonumber
    \boldsymbol{v}^T \boldsymbol{A}_{11} \boldsymbol{v} &\ge \sum_{K \in \mathcal{T}_h} \left( - \frac{1}{4\tau_K} \|v_h\|^2_{L^2(K)}
    + \varepsilon (\rho_{min}(\mathbb{D}_K) + \tau_K \rho_{min}(\mathbb{K}_K)) \|\nabla v_h\|^2_{L^2(K)} \right). 
\end{align}

Applying inverse inequality
\begin{equation}
    \|v_h\|_{L^2(K)} \le C_{inv} h_K \|\nabla v_h\|^2_{L^2(K)},
\end{equation}
we have
\begin{equation}
    \boldsymbol{v}^T \boldsymbol{A}_{11} \boldsymbol{v} \ge \sum_K \left( \varepsilon (\rho_{min}(\mathbb{D}_K) + \tau_K \rho_{min}(\mathbb{K}_K)) 
    - C^2_{inv} \frac{h_K^2}{4\tau_K} \right) \|\nabla v_h\|^2_{L^2(K)} > 0,
\end{equation}
provided the condition \eqref{cond-stab-2} holds.
Thus, $\boldsymbol{v}^T \boldsymbol{B} \boldsymbol{v} > 0$ for all $\boldsymbol{v} \ne \boldsymbol{0}$. So $\boldsymbol{B}$ is positive definite. 

Since $\boldsymbol{B}$ is a $Z$-matrix that is positive definite, it follows that $\boldsymbol{B}$ is a nonsingular $M$-matrix. 
Hence $\boldsymbol{B}^{-1} \ge 0$. Then,
\begin{align*}
    \tilde{\boldsymbol A}^{-1} = \begin{bmatrix}
    \boldsymbol{B}^{-1} & -\boldsymbol{B}^{-1} \boldsymbol{C}\\
    0 & \boldsymbol{I}
    \end{bmatrix} \ge 0,
\end{align*}
since $\boldsymbol{C} \le 0$.
Therefore, $\tilde{\boldsymbol A}$ is an $M$-matrix.

\end{proof}

Under the same assumptions in Lemma \ref{lemma-diagonal}, we have
\[M_{ii} \ge \sum_{K \in \omega_i} \left(\frac{2|K|}{(d+1)(d+2)} - \tau_K \frac{|K|}{(d+1)h_i^K}\|b\|_{\infty,K} \right),\]
and
\begin{align}
    \nonumber
    A_{ii}&= \sum_{K\in\omega_i}\int_K \left( \varphi_i{\boldsymbol b}\cdot \nabla\varphi_i+\tau_K({\boldsymbol b\cdot \nabla\varphi_i})^2  \right) \dif {\boldsymbol x} + |K| (\nabla\varphi_i)^T \varepsilon (\tau_K \mathbb{K}_K + \mathbb{D}_K) \nabla\varphi_i  \\
    & \le \sum_{K\in\omega_i} |K| \left( \frac{\|b\|_{\infty,K}}{(d+1) h_i^K} + \frac{\|b\|^2_{\infty,K}}{(h_i^K)^2} + \varepsilon \frac{\rho_{min}(\mathbb{D}_K)+\tau_K \rho_{min}(\mathbb{K}_K)}{(h_i^K)^2} \right).
\end{align}
Then we have the following lemma.
\begin{lemma}
Assume
\begin{equation}
    M_{ii}^{low} := \sum_{K \in \omega_i} \left(\frac{2|K|}{(d+1)(d+2)} - \tau_K \frac{|K|}{(d+1)h_i^K}\|b\|_{\infty,K} \right) > 0,
\end{equation}
and let
\begin{equation}
    A_{ii}^{up} := \sum_{K\in\omega_i} |K| \left( \frac{\|b\|_{\infty,K}}{(d+1) h_i^K} + \frac{\|b\|^2_{\infty,K}}{(h_i^K)^2} + \varepsilon \frac{\rho_{min}(\mathbb{D}_K)+\tau_K \rho_{min}(\mathbb{K}_K)}{(h_i^K)^2} \right).
\end{equation}
If 
\begin{equation}
    \label{temporalupperbound}
    \Delta t \le \frac{1}{1-\theta} \min_{i} \frac{M_{ii}^{low}}{A_{ii}^{up}},
\end{equation}
then 
\begin{equation}
    M_{ii} - \Delta t (1-\theta) A_{ii} \ge 0 \quad \text{ for all } i.
\end{equation}
\end{lemma}
\begin{proof}
Since $M_{ii} \ge M_{ii}^{low}$ and $A_{ii} \le A_{ii}^{up}$, we have
\begin{equation}
    M_{ii} - \Delta t (1-\theta) A_{ii} \ge M_{ii}^{low} - \Delta t (1-\theta) A_{ii}^{up} \ge 0
\end{equation}
whenever \eqref{temporalupperbound} holds. 
\end{proof}

\subsection{Metric tensor for DMP}
\label{subsec-metric}
We now develop the metric tensor to control the mesh so that the numerical solution satisfies DMP. Assume the temporal estimates from the previous section. As mentioned in the proof of Theorem 3.1 of \cite{li2013maximum}, $\tilde{\boldsymbol A}^{-1}\geq 0$ and $\tilde{\boldsymbol A}$ is an $M$-matrix. To see that $\tilde{\boldsymbol A}$ has nonnegative row sums, note that $\sum_j\varphi_j=1$. Therefore, the solution to (\ref{ODEsystem}) satisfies 
\begin{align}
    \underset{i=1,\dots,(N+1)N_v}{\max}(u_h)_i = \max \left\{0,\underset{i\in S({\boldsymbol f}^+)}{\max}(u_h)_i\right\},
\end{align}
where $S({\boldsymbol f}^+)$ is the set of indices where $f_i$ is nonnegative. Piecewise linear functions attain their maximum value at the vertices of element $K$. Hence, when $f({\boldsymbol x},t)\leq 0$, we have 
\begin{align}
    \underset{n = 0,\dots, N}{\max}\underset{{\boldsymbol x}\in\overline{\Omega}}{\max}\ u_h^n({\boldsymbol x}) = \max\left\{0,\underset{n = 0,\dots, N}{\max}\underset{{\boldsymbol x}\in\overline{\Omega}}{\max}\ u_h^n({\boldsymbol x}), \underset{{\boldsymbol x\in\overline{\Omega}}}{\max}\ u_h^0({\boldsymbol x})\right\}, 
\end{align}
where 
\begin{align}
    u_h^n = \sum_{j=1}^{N_{vi}}u_j^n\varphi_j({\boldsymbol x})+\sum_{j=N_{vi}+1}^{N_v} u_j^n\varphi_j({\boldsymbol x}), \ \ n=0,\dots, N. 
\end{align}

\begin{theorem}
    \label{thm-DMP}
Consider the fully discrete MM-SUPG scheme \eqref{MmatrixFormula}-\eqref{Frightside}. 
Assume that, for each element $K \in \mathcal{T}_h$, the following conditions hold:
\begin{enumerate}
\item[(i)] The elementwise tensors $\mathbb{D}_K$ and $\mathbb{K}_K$ are symmetric positive definite.
\item[(ii)] The meshes associated with $\mathbb{D}_K$ and $\mathbb{K}_K$ satisfy the anisotropic nonobtuse angle condition (ANAC).
\item[(iii)] The stabilization parameter $\tau_K$ satisfies \eqref{miicondition}, \eqref{aiicondition}, and \eqref{delta_ij_cond}.
\item[(vi)] The time step $\Delta t$ satisfies \eqref{temporallowerbound} and \eqref{temporalupperbound}.
\end{enumerate}
Then the coefficient matrices satisfy:
\begin{itemize}
\item[(a)] $\boldsymbol{M} - \Delta t(1-\theta)\boldsymbol{A}$ has nonnegative entries,
\item[(b)] $\widetilde{\boldsymbol{A}} = \boldsymbol{M} + \theta \Delta t \boldsymbol{A}$ is an $M$-matrix.
\end{itemize}
Consequently, the fully discrete MM-SUPG scheme satisfies the discrete maximum principle (DMP), i.e., the numerical solution preserves the bounds of the initial and boundary data.
\end{theorem}
\begin{proof}
    The results follows immediately from the previous lemmas. 
\end{proof}

\begin{corollary}
    Under the assumptions of Theorem \ref{thm-DMP}, the discrete maximum principle (DMP) holds if the following practical conditions are satisfied:
\begin{enumerate}
\item[(i)] \textbf{Mesh condition:}  
The mesh is sufficiently aligned with the anisotropy of the problem in the sense that it satisfies the anisotropic nonobtuse angle condition (ANAC) with respect to both $\mathbb{D}_K$ and $\mathbb{K}_K$.

\item[(ii)] \textbf{Stabilization condition:}  
The stabilization parameter $\tau_K$ is chosen such that
\[
\tau_K \|b\|_{\infty,K} \lesssim h_K,
\]
where $h_K$ denotes a representative element size.

\item[(iii)] \textbf{Balance condition:}  
The combined diffusion and stabilization effects dominate the convection contribution, i.e.,
\[
\rho_{\min}(\varepsilon \mathbb{D}_K) + \tau_K  \rho_{\min}(\varepsilon \mathbb{K}_K)
\gtrsim
h_K \|b\|_{\infty,K} + \tau_K \|b\|_{\infty,K}^2.
\]

\item[(iv)] \textbf{Time step condition:}  
The time step $\Delta t$ is chosen within a suitable range such that
\[
\Delta t \gtrsim \frac{1}{\theta} \frac{M_{ij}}{|A_{ij}|}
\quad \text{and} \quad
\Delta t \lesssim \frac{1}{1-\theta} \frac{M_{ii}}{A_{ii}},
\]
ensuring that $\widetilde{A} = M + \theta \Delta t A$ is an $M$-matrix and that $M - \Delta t(1-\theta)A$ remains nonnegative.
\end{enumerate}

Under these conditions, the fully discrete MM-SUPG scheme satisfies the discrete maximum principle.
\end{corollary}

\begin{remark}
    The corollary highlights that the DMP is achieved through a balance between mesh alignment, stabilization, and time stepping. In particular, the interaction between the mesh metric and the SUPG stabilization plays a key role in controlling oscillations and preserving monotonicity.
\end{remark}

\begin{remark}
    Since we use implicit Euler method in time for all computations (that is, $\theta =1$), the lower bound \eqref{temporallowerbound} for time step $\Delta t$ is practically removed. 
\end{remark}

Since two different metric tensors will be used on any given node, we follow  \cite{alauzet2011extension} and  take into consideration for every time step, we use metric intersection in order to generate the anisotropic mesh. 
\begin{figure}[thb]
    \centering
    \scalebox{0.75}{
    \includegraphics{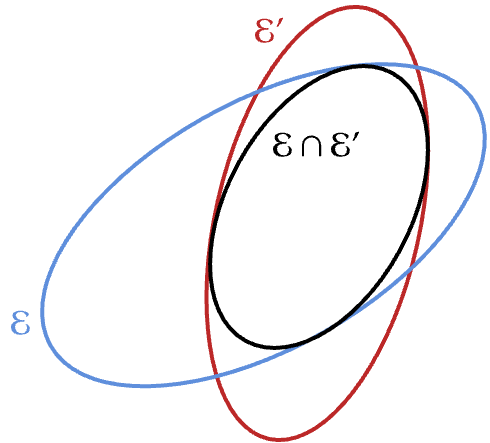}
    }
    \caption{Geometric illustration of intersecting ellipses $\mathcal{E}$ and $\mathcal{E}'$, resulting in $\mathcal{E}\cap \mathcal{E}'$.}
        \label{ellilpseIntersection}
\end{figure}
This is equivalent to finding the largest ellipsoid included in the intersection of the two ellipses $\mathcal{E}$ and $\mathcal{E}'$, the geometric representations of symmetric, positive definite $\mathcal{M}$ and $\mathcal{M}'$, respectively (see Figure \ref{ellilpseIntersection}). By our assumptions on $\mathcal{M}$ and $\mathcal{M}'$, $\mathcal{N} = \mathcal{M}^{-1}\mathcal{M}'$ is diagonalizable by the similarity transform. Denote the shared normalized eigenvectors of $\mathcal{M}$ and $\mathcal{M}'$ by ${\boldsymbol \lambda}_1$ and ${\boldsymbol \lambda}_2$ and write 
\begin{align}
    \label{M-int}
    \mathcal{M}_{int} = \mathcal{M} \cap \mathcal{M}' = \Lambda^{-T}\textrm{diag}\left(\max\{\alpha_1,\beta_1\},\max\{\alpha_2,\beta_2\}\right)\Lambda^{-1}.
\end{align}
Figure \ref{cooridnatetransformation} shows the coordinate transformation of a typical element $K$ to $\tilde{K}$ and $\check{K}$.
\begin{figure}[thb]
    \centering
    \tikzset{every picture/.style={line width=0.75pt}} %set default line width to 0.75pt        
    \begin{tikzpicture}[x=0.75pt,y=0.75pt,yscale=-1,xscale=1]
        %uncomment if require: \path (0,385); %set diagram left start at 0, and has height of 385

        %Flowchart: Merge [id:dp8796123094439892] 
        \draw   (222.07,97.85) -- (218,297.07) -- (138.96,195.8) -- cycle ;
        %Flowchart: Merge [id:dp57103384758654] 
        \draw   (489.15,190.09) -- (357.59,339.74) -- (362.47,211.37) -- cycle ;
        %Flowchart: Merge [id:dp8421802757428951] 
        \draw   (359.07,0.46) -- (474.01,163.22) -- (350.3,128.62) -- cycle ;
        %Straight Lines [id:da8543510950563746] 
        \draw    (233.33,170.47) -- (327.43,140.58) ;
        \draw [shift={(329.33,139.97)}, rotate = 522.37] [color={rgb, 255:red, 0; green, 0; blue, 0 }  ][line width=0.75]    (10.93,-3.29) .. controls (6.95,-1.4) and (3.31,-0.3) .. (0,0) .. controls (3.31,0.3) and (6.95,1.4) .. (10.93,3.29)   ;
        %Straight Lines [id:da8357421842931045] 
        \draw    (233.33,220.67) -- (327.37,240.1) ;
        \draw [shift={(329.33,240.5)}, rotate = 191.67] [color={rgb, 255:red, 0; green, 0; blue, 0 }  ][line width=0.75]    (10.93,-3.29) .. controls (6.95,-1.4) and (3.31,-0.3) .. (0,0) .. controls (3.31,0.3) and (6.95,1.4) .. (10.93,3.29)   ;

        % Text Node
        \draw (183.67,188.2) node [anchor=north west][inner sep=0.75pt]    {$K$};
        % Text Node
        \draw (376.92,94.13) node [anchor=north west][inner sep=0.75pt]    {$\tilde{K}$};
        % Text Node
        \draw (386.42,227.54) node [anchor=north west][inner sep=0.75pt]    {$\check{K}$};
        % Text Node
        \draw (266.33,118) node [anchor=north west][inner sep=0.75pt]    {$\mathbb{D}_K^{-1/2}{\boldsymbol x}$};
        % Text Node
        \draw (275.33,202.2) node [anchor=north west][inner sep=0.75pt]    {$\mathbb{K}_K^{-1/2}{\boldsymbol x}$};

    \end{tikzpicture}
    \caption{Coordinate transformation of a typical element $K$ to $\tilde{K}$ and $\check{K}$.}
    \label{cooridnatetransformation}
\end{figure}

\subsection{Numerical Results}

In this section we present our numerical findings for using different metric tensors. 
We start with the metric tensors $\mathcal{M}_{DMP} = c_1 \mathbb{D}^{-1}$ and $\mathcal{M}_{\mathbb{K}} = c_2 \mathbb{K}^{-1}$, where $c_1$ and $c_2$ are arbitray constants. 
Define the metric tensor 
\begin{equation}
    \mathcal{M}_{int}=\mathcal{M}_{DMP} \cap \mathcal{M}_{\mathbb{K}}.
\end{equation}
We use the implicit Euler method in time for all computations (that is, $\theta =1$) and compare the results with the classic Galerkin method with SUPG using fixed mesh (FM-SUPG).

%%%%%%%%%%%%%%%%%%%%%%          Example 3       %%%%%%%%%%%%%%%%%%%%%%%%%%%%%%%%%%%%%%%%

\begin{exam}
\label{example3}
Consider the manufactured solution on the domain $\Omega_T = [0,1]^2\times[0, 0.3]$ to (\ref{modeleqn}) given by 
\begin{align}
    u(x,y,t) = (1+\exp(C(x+y-t))^{-1}, \ \ C>0,
\end{align}
with diffusion tensor defined by 
\begin{align}
    \label{diffusionmatrix1}
    \mathbb{D} = \begin{pmatrix}
    50.5 & 49.5\\
    49.5 & 50.5
    \end{pmatrix},
\end{align}
and time-dependent flow given by 
\begin{align}
    \label{ex3-flow}
    {\boldsymbol b}(x,y,t) = (x+t,y+t). 
\end{align}
In the computations, we have used $\varepsilon = 10^{-6}$ and $\Delta t = 10^{-2}$ for this example.
Based on the maximum principle, the solution should satisfy $\underset{\Omega_T}{\min} \ u = 0$.

The diffusion tensor has eigenvalues $\lambda_1 = 100$ and $\lambda_2 = 1$, 
which indicates that the diffusion in horizontal direction is 100 times faster than that in the vertical direction. 
For the given flow in \eqref{ex3-flow}, the tensor $\mathbb{K}$ defined in \eqref{Kmatrix} is the same as $\mathbb{D}$.
Therefore, the metric tensor $\mathcal{M}_{DMP}$ is the same as $\mathcal{M}_{\mathbb{K}}$ for this example.
Hence, for this example, $\mathcal{M}_{int}$ is also the same as $\mathcal{M}_{DMP}$.

% Table \ref{ex3-table1} presents the minimum values of the numerical solutions obtained from different meshes.
% As can be seen, the results from $\mathcal{M}_{int}$ and $\mathcal{M}_{DMP}$ are nearly identical. 
% Even the results from the fixed mesh show satisfaction of DMP for the bottom shelf. 
Figure \ref{ex3-sideview} shows the side view of the numerical solutions using different metric tensors. 
As can be seen, the graph in Figure \ref{ex3-sideview}(a) shows significant artificial oscillation near the top shelf when using the fixed mesh approach. 
On the other hand, as shown in Figure \ref{ex3-sideview}(b), when using $\mathcal{M}_{int}$ or $\mathcal{M}_{DMP}$,
the artificial oscillation at the top shelf is essentially eliminated. 

% \begin{table}[thb]
%     \caption{Example \ref{example3}. $\underset{\Omega_T}{\min} \ u$ with $\Delta t = 10^{-2}$. }
%     \centering 
%     \begin{tabular}{|c|c|c|l|}
%         \hline
%         $N$ & $\mathcal{M}_{int}$ & $\mathcal{M}_{DMP}$ & Fixed mesh      \\ \hline
%         16200 & -4.4374e-05            & -4.4373e-05  & -2.8664e-07      \\ \hline
%         24200 & -3.1365e-06            & -3.1367e-06  & -2.2583e-08      \\ \hline
%         34322 & -1.5736e-07            & -1.5735e-07  & -1.1880e-09      \\ \hline
%     \end{tabular}
%     \label{ex3-table1}
% \end{table}

\begin{figure}[thb]
    \centering
    \hbox{
    \begin{minipage}[t]{3.0in}
    \includegraphics[width=3.0in]{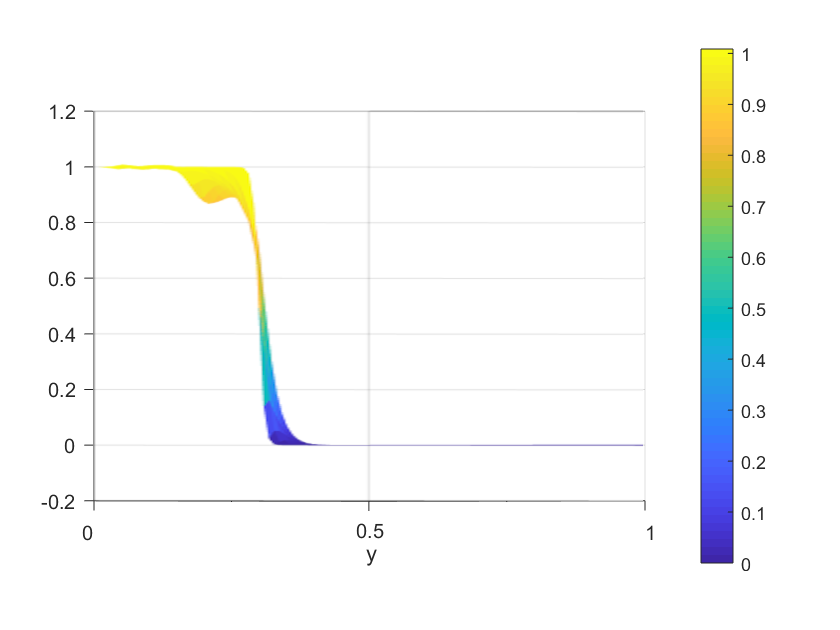}
    \centerline{(a): Fixed mesh}
    \end{minipage}
    \hspace{10mm}
    \begin{minipage}[t]{3.0in}
    \includegraphics[width=3.0in]{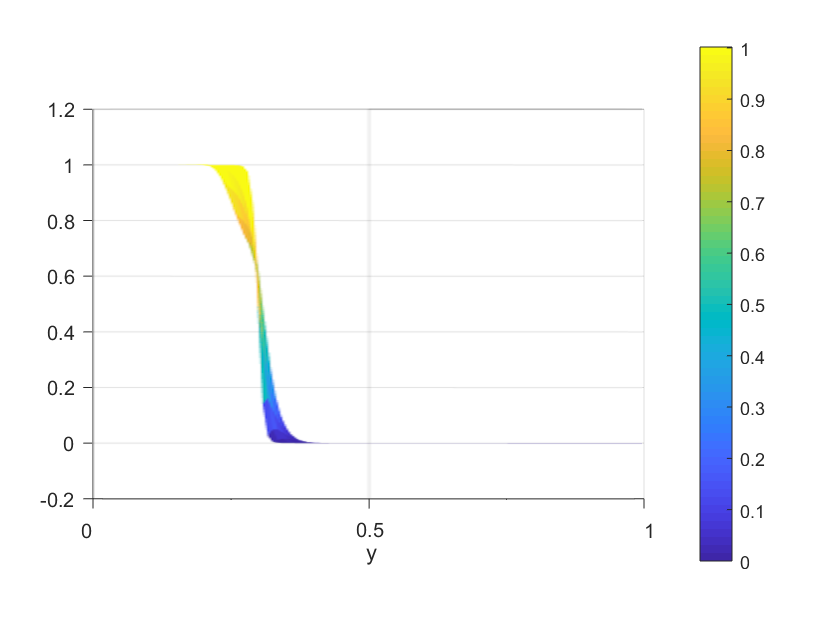}
    \centerline{(b): $\mathcal{M}_{int}$}
    \end{minipage}
    }
    \caption{Example \ref{example3} -- side view of solutions (scaled-up) for a fixed mesh and $\mathcal{M}_{int}$ with $N = 34322$, $\varepsilon = 10^{-6}$, and $\Delta t = 10^{-2}$.} 
    \label{ex3-sideview}
\end{figure}

\end{exam}

%%%%%%%%%%%%%%%%%%%%%%          Example 4       %%%%%%%%%%%%%%%%%%%%%%%%%%%%%%%%%%%%%%%%

\begin{exam}
\label{example4}
In this example, we consider the metric tensor $\mathcal{M}_{DMP+adapt}$ developed in \cite{li2010anisotropic} that is 
based on both DMP satisfaction as well as minimization of the $H^1$-seminorm of the interpolation error. 
Correspondingly, our intersection metric tensor is $\mathcal{M}_{int+adapt} = \mathcal{M}_{DMP+adap} \cap \mathcal{M}'$ where $\mathcal{M}' = \theta' \mathbb{K}^{-1}$. 
The anisotropic diffusion tensor is given by 
\begin{align}
    \mathbb{D} = \begin{pmatrix}
    1&1\\1&3
    \end{pmatrix},
\end{align}
and the flow is given by 
\begin{align}
    {\boldsymbol b}(x,y,t) = (x+y+t+8,x+3y+t+16),
\end{align}
for $\Omega_T = [0,1.5]^2\times[0,0.015]$. 
For the given ${\boldsymbol b}$ and $\mathbb{D}$, we have 
\begin{align}
    \mathbb{K} = \begin{pmatrix}
    2&4\\4&10
    \end{pmatrix}.
\end{align}

We first use $N=20,402$ elements in the mesh and compute the numerical solutions for $\varepsilon=10^{-3}$ with time step size $\Delta t = 10^{-4}$.
Figure \ref{ex4-sideview} shows the side view of the solution using fixed mesh, $\mathcal{M}_{DMP+adapt}$ mesh and $\mathcal{M}_{int+adapt}$ mesh, respectively. 
As can be seen, $\mathcal{M}_{DMP+adapt}$ mesh and $\mathcal{M}_{int+adapt}$ mesh provide comparable results and both perform much better than the fixed mesh SUPG approach. 

\begin{figure}[!thb]
    \centering
    \hbox{
    \begin{minipage}[t]{2.0in}
    \includegraphics[width=2.0in]{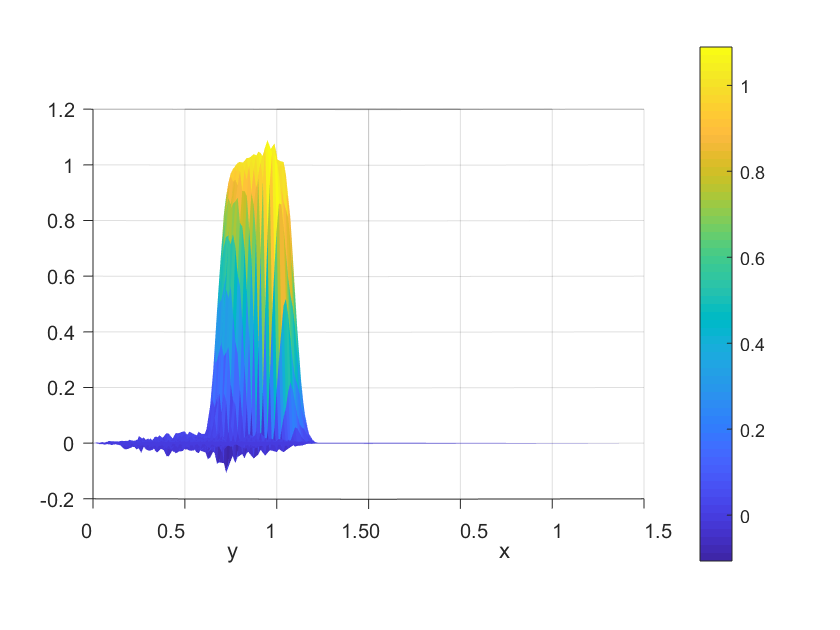}
    \centerline{(a): Fixed mesh}
    \end{minipage}
    \begin{minipage}[t]{2.0in}
    \includegraphics[width=2.0in]{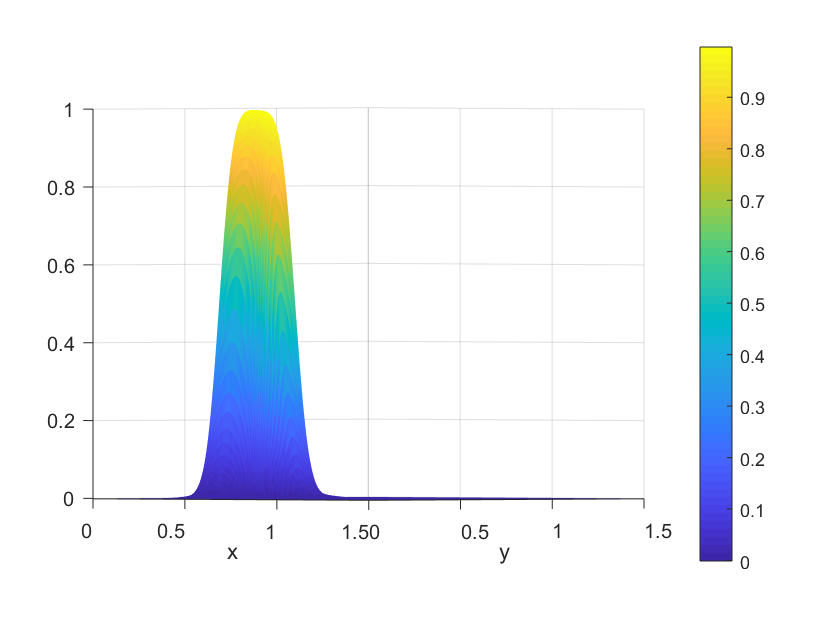}
    \centerline{(b): $\mathcal{M}_{DMP+adapt}$}
    \end{minipage}
    \begin{minipage}[t]{2.0in}
    \includegraphics[width=2.0in]{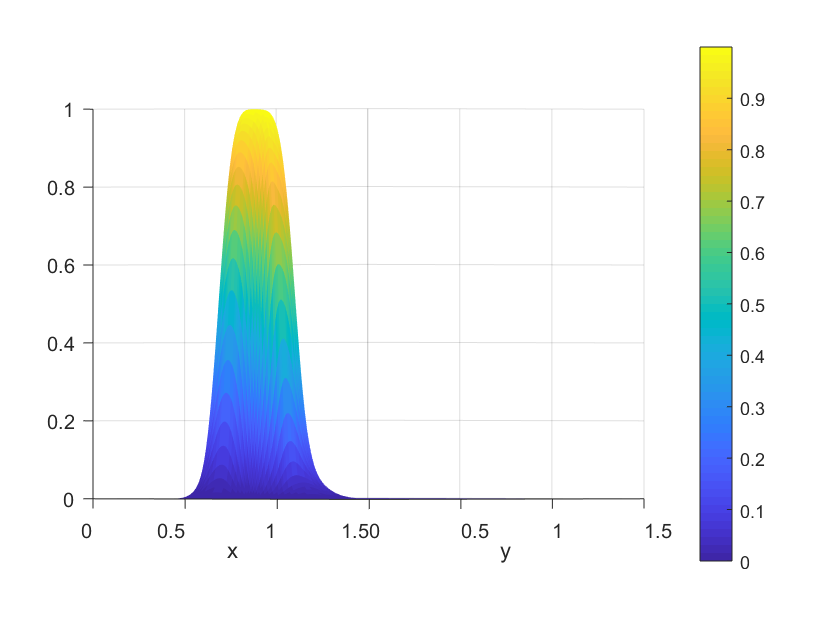}
    \centerline{(c): $\mathcal{M}_{int+adapt}$}
    \end{minipage}
    }
    \caption{Example \ref{example4} -- side view with $\varepsilon = 10^{-3}$ for different metric tensors and $\Delta t = 10^{-4}$, $N = 20,402$.}
    \label{ex4-sideview}
\end{figure}

To see the advantage of $\mathcal{M}_{int+adapt}$ mesh over $\mathcal{M}_{DMP+adapt}$ mesh, we compute the numerical solutions at different $N$ values with larger time step $\Delta t = 10^{-3}$.
The minimum values of the numerical solutions are summarized in Table \ref{ex4-table}. As can be seen, except on the course mesh, $\mathcal{M}_{int+adapt}$ provides the best results. 

\begin{table}
    \caption{Example \ref{example4}. Results for $\mathcal{M}_{int+adapt}$, $\mathcal{M}_{DMP+adapt}$, and fixed mesh for varying $N$ with $\varepsilon = 10^{-3}$ and $\Delta t = 10^{-3}$.}
    \centering
    \begin{tabular}{|c|c|c|c|}
    \hline
    $N$ & Fixed Mesh & $\mathcal{M}_{DMP+adapt}$  & $\mathcal{M}_{int+adapt}$  \\ \hline
    3362  & -2.4050e-01             & -2.3330e-03             & -3.6297e-03      \\ \hline
    13122 & -1.6762e-01             & -1.3796e-03             & -6.8634e-04      \\ \hline
    29282 & -1.0677e-01             & -9.3509e-04             & -1.8547e-04      \\ \hline
    \end{tabular}
    \label{ex4-table}
\end{table}
    
\end{exam}

\section{Summary and Conclusions}

In this paper, we have developed a moving mesh streamline upwind Petrov-Galerkin (MM-SUPG) method that combines the moving mesh partial differential equation (MMPDE) approach with the SUPG stabilization for convection-diffusion problems. 
The key feature of the proposed method is not only the combination of these two techniques, but also their interaction: the evolving mesh modifies the local element geometry and hence influences both the direction and magnitude of the SUPG stabilization, 
while the stabilization improves robustness in convection-dominated regimes.

The MM-SUPG method leverages the complementary roles of the two components. The moving mesh method concentrates mesh elements in regions with sharp gradients, improving resolution and computational efficiency, 
while the SUPG method stabilizes the numerical solution by adding streamline diffusion. The method is applicable to both time-dependent and time-independent problems, and is designed to handle both isotropic and anisotropic diffusion tensors.

For the isotropic diffusion case, the MM-SUPG method exhibits approximately $O(h^{1})$ convergence in the $H^1$-seminorm and $O(h^{1.5})$ convergence in the $L^2$ norm (see Figures \ref{ex2-timeindependent} and \ref{ex2-timedependent}). 
Numerical experiments show that, for convection-dominated problems with sharp layers, moving mesh with the standard Galerkin method (MM-FEM) may still exhibit noticeable oscillations on moderately refined meshes (see Figure \ref{ex1-sideview}). 
The use of SUPG stabilization within the moving mesh framework reduces these oscillations and leads to improved stability. Comparisons in Examples \ref{example1} and \ref{example2} indicate that MM-SUPG improves robustness over fixed mesh SUPG (FM-SUPG), 
and provides comparable accuracy to MM-FEM in cases where mesh adaptation alone already captures the solution well.

For the anisotropic diffusion case, we have derived sufficient conditions on the mesh geometry and time step size under which the fully discrete scheme satisfies the discrete maximum principle (DMP). 
The analysis is based on quantitative bounds for the convection terms rather than sign arguments, and therefore avoids restrictive assumptions on the divergence of the velocity field. 
Based on the diffusivity tensor $\mathbb{D}$ and the flow vector $\boldsymbol{b}$, we introduce a weighted tensor $\mathbb{K}$ in \eqref{Kmatrix}, which serves as an analytical tool for incorporating the interaction between convection and diffusion. 
The tensor $\mathbb{K}$ is not guaranteed to be symmetric positive definite in general; thus, the DMP analysis is established under structural assumptions on $\boldsymbol{b}$ and $\mathbb{D}$ such that its elementwise form is of SPD type.

To further improve monotonicity properties, we construct a metric tensor $\mathcal{M}_{int}$ in \eqref{M-int} through metric intersection, combining information from both $\mathbb{D}$ and $\mathbb{K}$. 
Numerical results in Examples \ref{example3} and \ref{example4} show that meshes generated using $\mathcal{M}_{int}$ improve DMP behavior, as evidenced by reduced undershoots, while differences in standard error measures remain modest.

In summary, the MM-SUPG method provides a stable and effective approach for convection-dominated convection-diffusion problems by combining mesh adaptation and stabilization in a unified framework. 
Under suitable conditions on the mesh, problem data, and time step size, the method can produce solutions that satisfy the discrete maximum principle, particularly in anisotropic diffusion settings. 
The results demonstrate that the interaction between mesh adaptation and stabilization plays a key role in achieving both accuracy and improved monotonicity properties.

% \bibliographystyle{ieeetr}
% \bibliography{ref.bib}
%\bibliographystyle{unsrt}

\end{document}